\documentclass{iopart}
\usepackage{iopams}
\eqnobysec
\usepackage{graphicx}
\usepackage{epstopdf}
\newtheorem{lemma}{Lemma}[section]
\newtheorem{theorem}[lemma]{Theorem}
\newtheorem{remark}[lemma]{Remark}

\newtheorem{proposition}[lemma]{Proposition}

\newtheorem{algorithm}[lemma]{Algorithm}

\newtheorem{Assumption}[lemma]{Assumption}
\newcommand{\proof} [1]
   { {\bf Proof.} #1 \hfill\opensquare \\}
\renewcommand{\theequation}{\mbox{\arabic{section}.\arabic{equation}}}
\def\p{\partial}
\def\d{\delta}

\def\l{\langle}
\def\r{\rangle}
\def\N{\mathcal N}
\def\R{\mathcal R}
\def\D{\mathcal D}
\def\X{\mathcal X}
\def\Y{\mathcal Y}

\def\d{\delta}
\def\la{\lambda}
\def\X{\mathcal X}
\def\Y{\mathcal Y}

\begin{document}
\jl{5}
\title[]
{Landweber iteration of Kaczmarz type with general non-smooth convex penalty functionals}

\author{Qinian Jin$^1$ and Wei Wang$^{2,3}$}

\medskip

\address{$^1$Mathematical Sciences Institute, Australian National University,
Canberra, ACT 0200, Australia}
\address{$^2$College of Mathematics, Physics and Information Engineering, Jiaxing University, Zhejiang 314001, China}
\address{$^3$School of Mathematical Sciences, Fudan University, Shanghai 200433, China}

%\medskip

\eads{Qinian.Jin@anu.edu.au and weiwangmath@gmail.com }

%%%%%%%%%%%%%%%%%%%%%%%%%%%%%%%%%%%%%%%%%%%%%%%%%%%%%%%%%%%%%%%%%%%%%%%%%%%%%%%%

%\keywords{Ill-posed problems, inverse problems, regularization,
%Hilbert scales, implicit iteration method, order optimal error
%bounds, general smoothness conditions, operator monotone functions}

\begin{abstract}
   The determination of solutions of many inverse problems usually requires a set of measurements
   which leads to solving systems of ill-posed equations. In this paper we propose the Landweber iteration
   of Kaczmarz type with general uniformly convex penalty functional. The method is formulated by using tools
   from convex analysis. The penalty term is allowed to be non-smooth to include the $L^1$ and total variation (TV)
   like penalty functionals, which are significant in reconstructing special features of solutions such as
   sparsity and piecewise constancy in practical applications. Under reasonable conditions, we establish the convergence
   of the method. Finally we present numerical simulations on tomography problems and parameter identification
   in partial differential equations to indicate the performance.
\end{abstract}

%\ams{65J15, 65J20, 47H17}

%\submitto{\IP}

%\maketitle

%
\oddsidemargin 15mm
\evensidemargin 15mm
\topmargin  0 pt
%\textheight=8.6 true in
%\textwidth=5.2 true in

%\vskip 0.5 cm

%%%%%%%%%%%%%%%%%%%%%% -- SECTION 1 -- %%%%%%%%%%%%%%%%%%%%%%%%%%%%%%%%%%%%%%%%
%%%%%%%%%%%%%%%%%%%%%% -- SECTION 1 -- %%%%%%%%%%%%%%%%%%%%%%%%%%%%%%%%%%%%%%%%
%\section{Introduction}\label{s1}
%

\section{Introduction}

Landweber iteration is one of the most well-known regularization methods for solving inverse problems formulated
in Hilbert spaces. A complete account on this method for linear inverse problems can be found in \cite{EHN96}
including the convergence analysis and its various accelerated versions. A nonlinear
version of Landweber iteration was proposed in \cite{HNS95} for solving nonlinear inverse problems, where an elegant
convergence analysis was present. Although Landweber iteration converges slowly, it still receives a lot of attention
because it is simple to implement and is robust with respect to noise.

The classical Landweber iteration in Hilbert spaces, however, has the tendency to over-smooth solutions which
makes it difficult to capture the special features of the sought solutions such as sparsity and piecewise constancy.
It is therefore necessary to reformulate this method either in Banach space setting or in a manner that modern non-smooth
penalty functionals, such as the $L^1$ and total variation like functionals, can be incorporated.

Let $A: \X \to \Y$ be a linear compact operator between two Banach spaces $\X$ and $\Y$ with norms $\|\cdot\|$ whose dual spaces are denoted
by $\X^*$ and $\Y^*$ respectively. Some recent advances on Landweber iteration for linear inverse problems
\begin{equation}\label{Linear}
A x= y
\end{equation}
in Banach space setting have been reported using only the noisy data $y^\d$ satisfying
\begin{equation*}
\|y^\d-y\|\le \d
\end{equation*}
with a small known noise level $\d>0$.  In particular, when $\X$ is uniformly smooth and uniformly convex, by virtue of the duality mappings,
a version of Landweber iteration for solving (\ref{Linear}) was proposed in \cite{SLS2006}. Although the method
excludes the use of the $L^1$ and total variation like penalty functionals, new ideas were introduced
in \cite{SLS2006} which promote the study of Landweber iteration in modern setup. Recently a version of Landweber iteration
was proposed in \cite{BH2012} using non-smooth uniformly convex penalty functionals. Let $\Theta:\X\to (-\infty, \infty]$ be a proper,
lower semi-continuous, uniformly convex functional, then the method in \cite{BH2012} reads as
\begin{eqnarray}\label{Land2}
\left\{\begin{array}{lll}
\xi_{n+1} =\xi_n -\mu_n A^* J_r(A x_n-y^\d),\\
x_{n+1} =\arg\min_{x\in \X} \left\{\Theta(x) -\l \xi_{n+1}, x\r\right\},
%x_{n+1} =\nabla \Theta^*(\xi_{n+1}),
\end{array}\right.
\end{eqnarray}
where $A^*:\Y^*\to \X^*$ denotes the adjoint of $A$, $J_r$ with $1<r<\infty$ is the duality mapping of $\Y$ with gauge
function $t\to t^{r-1}$, $\{\mu_n\}$ are suitable chosen step-lengths,
%\begin{equation*}
%x_{n+1} =\arg\min_{x\in \X} \left\{\Theta(x) -\l \xi_{n+1}, x\r\right\},
%\end{equation*}
and $\l\cdot, \cdot\r$ denotes the duality pairing between $\X^*$ and $\X$. The method (\ref{Land2}) reduces to the one in \cite{SLS2006}
when taking $\Theta(x):=\|x\|^p/p$ with $1<p<\infty$. However, (\ref{Land2})
has more freedom on $\Theta$ so that it can be used to detect special features of solutions.

The convergence analysis of (\ref{Land2}) is given in \cite{BH2012} when it is terminated by the discrepancy principle
\begin{equation}\label{DP1}
\|A x_{n_\d}-y^\d\| \le \tau \d <\|A x_n-y^\d\|, \qquad 0\le n<n_\d
\end{equation}
with $\tau>1$. The argument in \cite{BH2012}, however, requires that $\mbox{int}(D(\Theta))$, the interior of $D(\Theta)$, must be non-empty
and that (\ref{Linear}) must have a solution in $\mbox{int}(D(\Theta))$. These conditions are indeed quite restrictive; for instance,
the domain of the total variation like functional
\begin{equation*}
\Theta(x):=\frac{1}{2\beta} \int_\Omega |x(\omega)|^2 d\omega + \int_\Omega |D x|
\end{equation*}
with $\beta>0$ over the space $\X:=L^2(\Omega)$ on a bounded domain $\Omega\subset {\mathbb R}^d$ does not have any interior
point in $L^2(\Omega)$, where
\begin{equation*}
\int_\Omega |D x| :=\sup\left\{ \int_\Omega x \, \mbox{div} f d\omega:
f\in C_0^1(\Omega; {\mathbb R}^N) \mbox{ and }
\|f\|_{L^\infty(\Omega)}\le 1\right\}
\end{equation*}
denotes the total variation of $x$ over $\Omega$ (\cite{G84}). Therefore, the theoretical result in \cite{BH2012} can not be applied to this important
penalty functional.

It is natural to ask if the convergence of (\ref{Land2}) can be proved without assuming $\mbox{int}(D(\Theta))\ne \emptyset$.
An affirmative answer would theoretically justify the applicability of (\ref{Land2}) to a wider class of penalty functionals $\Theta$
including the total variation like functionals. The control of $\{\xi_n\}$ presents one of the major challenges. The analysis in \cite{BH2012} is based on proving
the boundedness of $\{\xi_n\}$ in $\X^*$ which consequently enforces to assume that $\mbox{int}(D(\Theta))\ne \emptyset$.
We observe that the boundedness of $\{\xi_n\}$ is not essential in the convergence analysis, the most essential ingredient
is to control $\l \xi_n, x_n-\hat x\r$ for any solution $\hat x$ of (\ref{Linear}). Due to the lack of monotonicity on the residual
$\|A x_n-y^\d\|$, it turns out to be difficult to consider $\l \xi_n, x_n-\hat x\r$ for all $n\ge 0$. Fortunately, with a careful
chosen subsequence $\{n_k\}$ of integers, we can derive what we expect on $\l \xi_{n_k}, x_{n_k}-\hat x\r$ which together with
some monotonicity results enables us to prove a stronger result, i.e. $x_{n_\d}$ converges to a solution of (\ref{Linear}) in Bregman distance.

Instead of considering (\ref{Land2}) for solving (\ref{Linear}) directly, we consider a more general setup in which
(\ref{Land2}) is extended for solving linear as well as nonlinear inverse problems. Instead of studying a single equation, we consider the system
\begin{equation}\label{sys0}
F_i(x) =y_i, \qquad i=0, \cdots, N-1
\end{equation}
consisting of $N$ equations, where, for each $i=0, \cdots, N-1$, $F_i: D(F_i)\subset \X\to \Y_i$ is an operator between two
Banach spaces $\X$ and $\Y_i$. Such systems arise naturally in many practical applications including various tomography techniques
using multiple exterior measurements. By introducing
\begin{equation*}
F := (F_0, \cdots, F_{N-1}) : \D=\bigcap_{i=0}^{N-1} D(F_i) \to \Y_0\times \cdots \times \Y_{N-1}
\end{equation*}
and
\begin{equation*}
y:=(y_0, \cdots, y_{N-1}),
\end{equation*}
the system (\ref{sys0}) could be reformulated as a single equation $F(x) = y$. One might consider extending (\ref{Land2}) to
solve $F(x)=y$ directly. This procedure, however, becomes inefficient if $N$ is large because it destroys the special structure of (\ref{sys0}) and
results in an equation requiring huge memory to save the intermediate computational results. Therefore, it seems advantageous
to use the Kaczmarz-type methods, which cyclically consider each equation in (\ref{sys0}) separately and hence require only reasonable
memory consumption.

Some Landweber-Kaczmarz methods were formulated in \cite{KS2002,HLS2007} for solving the system (\ref{sys0}) when $\X$ and $\Y_i$ are
Hilbert spaces,  and the numerical results indicate that artefacts can appear in the reconstructed solutions
due to oversmoothness. Recently a Landweber-Kaczmarz method was proposed in \cite{LA2012} for solving (\ref{sys0})
in Banach space setting in the spirit of \cite{SLS2006} and hence the possible use of the
$L^1$ and total variation like penalty functionals is excluded. Furthermore, the convergence analysis in \cite{LA2012} unfortunately contains an error
(see the first line on page 12 in \cite{LA2012}). In this paper, we propose a Landweber iteration of Kaczmarz type
in which (\ref{Land2}) is adapted to solve each equation in (\ref{sys0}) and thus general non-smooth uniformly convex penalty
functionals $\Theta$ are incorporated into the method with the hope of removing artefacts and of capturing special features of solutions.
We give the detailed convergence analysis of our method. It is worthy pointing out that our analysis does not require the interior
of $D(\Theta)$ be nonempty and therefore the convergence result applies for the total variation like penalty functionals.

This paper is organized as follows. In section 2 we give some preliminary results from convex analysis. In section 3, we first formulate
the Landweber iteration of Kaczmarz type with general uniformly convex penalty term for solving the system (\ref{sys0}), and then
present the detail convergence analysis. In section 4 we give the proof of an important proposition which plays an important role
in section 3. Finally, in section 5 we present some numerical simulations on tomography problems in imaging and parameter identification
in partial differential equations to test the performance of the method.

\section{Preliminaries}

Let $\X$ be a Banach space with norm $\|\cdot\|$. We use $\X^*$ to denote its dual space, and for any $x\in \X$
and $\xi \in \X^*$ we write $\l \xi, x\r=\xi (x)$ for the duality pairing. If $\Y$ is another Banach space
and $A: \X\to \Y$ is a bounded linear operator, we use $A^*: \Y^*\to \X^*$ to denote its adjoint, i.e.
$\l A^* \zeta, x\r=\l \zeta, A x\r$ for any $x\in \X$ and $\zeta\in \Y^*$. Let $\N(A)=\{x\in \X: A x=0\}$ be
the null space of $A$ and let
\begin{equation*}
\N(A)^\perp :=\{\xi\in \X^*: \l \xi, x\r =0 \mbox{ for all } x\in \N(A)\}
\end{equation*}
be the annihilator of $\N(A)$. When $\X$ is reflexive, there holds
\begin{equation*}
\N(A)^\perp =\overline{\R(A^*)},
\end{equation*}
where $\overline{\R(A^*)}$ denotes the closure of $\R(A^*)$, the range space of $A^*$, in $\X^*$.

Given a convex function $\Theta: \X \to (-\infty, \infty]$, we use
\begin{equation*}
D(\Theta):=\{x\in \X: \Theta(x)<+\infty\}
\end{equation*}
to denote its effective domain. It is called proper if $D(\Theta)\ne \emptyset$. The subgradient of $\Theta$ at
$x\in \X$ is defined as
\begin{equation*}
\p \Theta(x):=\{\xi\in \X^*: \Theta(z)-\Theta(x)-\l \xi, z-x\r \ge 0 \mbox{ for all } z\in \X\}.
\end{equation*}
The multi-valued mapping $\p \Theta: \X\to 2^{\X^*}$ is called the subdifferential of $\Theta$. We set
\begin{equation*}
D(\p\Theta):=\{x\in D(\Theta): \p \Theta(x)\ne \emptyset\}.
\end{equation*}
For $x \in D(\p \Theta)$ and $\xi\in \p \Theta(x)$ we define (\cite{Br1967})
\begin{equation*}
D_\xi \Theta(z,x):=\Theta(z)-\Theta(x)-\l \xi, z-x\r, \qquad \forall z\in \X
\end{equation*}
which is called the Bregman distance induced by $\Theta$ at $x$ in the direction $\xi$. Clearly $D_\xi \Theta(z,x)\ge 0$
and
\begin{equation}\label{4.3.1}
D_\xi \Theta(x_2,x)-D_\xi \Theta(x_1, x) =D_{\xi_1} \Theta(x_2,x_1) +\l \xi_1-\xi, x_2-x_1\r
\end{equation}
for all $x, x_1\in D(\p \Theta)$, $\xi\in \p \Theta(x)$, $\xi_1\in \p \Theta(x_1)$ and $x_2\in \X$.

Bregman distance can be used to obtain information under the Banach space norm when $\Theta$ has stronger convexity.
A proper convex function $\Theta: \X \to (-\infty, \infty]$ is called uniformly convex if there is a continuous increasing
function $h:[0, \infty) \to [0, \infty)$, with the property that $h(t)=0$ implies $t=0$, such that
\begin{equation*}
\Theta(\la \bar x +(1-\la) x) +\la (1-\la) h(\|\bar x-x\|) \le \la \Theta(\bar x) +(1-\la) \Theta(x)
\end{equation*}
for all $\bar x, x\in \X$ and $\la\in (0,1)$. If $h$ can be taken as $h(t)=c_0 t^p$ for some $c_0>0$ and $p\ge 2$, then $\Theta$
is called $p$-convex. It can be shown that if $\Theta$ is uniformly convex then
\begin{equation*}
D_\xi \Theta(\bar x,x) \ge h(\|\bar x-x\|)
\end{equation*}
for all $\bar x \in \X$, $x\in D(\p \Theta)$ and $\xi\in \p \Theta(x)$. In particular, if $\Theta$ is $p$-convex with
$h(t) =c_0 t^p$, then
\begin{equation}\label{pconv}
D_\xi \Theta(\bar x,x) \ge c_0 \|\bar x-x\|^p
\end{equation}
for all $\bar x\in \X$, $x\in D(\p \Theta)$ and $\xi\in \p \Theta(x)$.

For a proper, lower semi-continuous, convex function $\Theta: \X \to (-\infty, \infty]$, its Legendre-Fenchel conjugate
is defined by
\begin{equation*}
\Theta^*(\xi):=\sup_{x\in \X} \left\{\l\xi, x\r -\Theta(x)\right\}, \quad \xi\in \X^*.
\end{equation*}
It is well known that $\Theta^*$ is also proper, lower semi-continuous, and convex. If, in addition,
$\X$ is reflexive, then
\begin{equation}\label{3.18.1}
\xi\in \p \Theta(x) \Longleftrightarrow x\in \p \Theta^*(\xi) \Longleftrightarrow \Theta(x) +\Theta^*(\xi) =\l \xi, x\r.
\end{equation}
When $\Theta$ is $p$-convex satisfying (\ref{pconv}) with $p\ge 2$, it follows from \cite[Corollary 3.5.11]{Z2002}
that $D(\Theta^*)=\X^*$, $\Theta^*$ is Fr\'{e}chet differentiable and its gradient $\nabla \Theta^*: \X^*\to \X$ satisfies
\begin{equation}\label{3.29.1}
\|\nabla \Theta^*(\xi_1)-\nabla \Theta^*(\xi_2) \|\le \left(\frac{\|\xi_1-\xi_2\|}{2c_0}\right)^{\frac{1}{p-1}},
\quad \forall \xi_1, \xi_2\in \X^*.
\end{equation}
Moreover
\begin{equation}\label{6.21.1}
\Theta^*(\xi_2)-\Theta^*(\xi_1) -\l \xi_2-\xi_1, \nabla \Theta^*(\xi_1)\r \le \frac{1}{p^*(2c_0)^{p^*-1}} \|\xi_2-\xi_1\|^{p^*}
\end{equation}
for any $\xi_1, \xi_2\in \X^*$, where $p^*$ is the number conjugate to $p$, i.e. $1/p+1/p^*=1$.
By the subdifferential calculus, there also holds
\begin{equation}\label{8.21.1}
x=\nabla \Theta^*(\xi) \Longleftrightarrow x =\arg \min_{z\in \X} \left\{ \Theta(z) -\l \xi, z\r\right\}.
\end{equation}

On a Banach space $\X$, we consider for $1<r<\infty$ the convex function $x\to \|x\|^r/r$.
Its subdifferential at $x$ is given by
\begin{equation*}
J_r^\X(x):=\left\{\xi\in \X^*: \|\xi\|=\|x\|^{r-1} \mbox{ and } \l \xi, x\r=\|x\|^r\right\}
\end{equation*}
which gives the duality mapping $J_r^\X: \X \to 2^{\X^*}$ of $\X$ with gauge function $t\to t^{r-1}$.
The duality mapping $J_r^\X$, for each $1<r<\infty$, is single valued and uniformly continuous on bounded sets if
$\X$ is uniformly smooth in the sense that its modulus of smoothness
\begin{equation*}
\rho_{\X}(s) := \sup\{\|\bar x+ x\|+\|\bar x-x\|- 2 : \|\bar x\| = 1,  \|x\|\le s\}
\end{equation*}
satisfies $\lim_{s\searrow 0} \frac{\rho_{\X}(s)}{s} =0$.

In many practical applications, proper, weakly lower semi-continuous, $p$-convex functions can be easily
constructed. For instance, consider $\X=L^p(\Omega)$, where $2\le p<\infty$
and $\Omega$ is a bounded domain in ${\mathbb R}^d$. It is known that the functional
\begin{equation*}
\Theta_0(x) := \int_\Omega |x(\omega)|^p d\omega
\end{equation*}
is $p$-convex on $L^p(\Omega)$. We can construct the new $p$-convex functions
\begin{equation}\label{eq:7.27}
\Theta(x):=\mu \int_\Omega |x(\omega)|^p d\omega + a \int_\Omega |x(\omega)| d\omega
+b \int_\Omega |D x|,
\end{equation}
where $\mu>0$, $a, b\ge 0$, and $\int_\Omega|D x|$ denotes the total variation of $x$ over $\Omega$.
For $a=1$ and $b=0$ the corresponding function is useful for sparsity reconstruction (\cite{T96});
while for $a=0$ and $b=1$ the corresponding function is useful for detecting the discontinuities,
in particular, when the solutions are piecewise-constant (\cite{ROF92}).

\section{Landweber iteration of Kaczmarz type}\label{Landweber}
\setcounter{equation}{0}

We consider the system (\ref{sys0}), i.e.
\begin{equation}\label{sys}
F_i(x) =y_i, \qquad i=0, \cdots, N-1
\end{equation}
consisting of $N$ equations, where, for each $i=0, \cdots, N-1$, $F_i: D(F_i)\subset \X\to \Y_i$ is an
operator between two Banach spaces $\X$ and $\Y_i$. We will assume that
\begin{equation*}
\D:=\bigcap_{i=0}^{N-1} D(F_i)\ne \emptyset
\end{equation*}
and each $F_i$ is Fr\'{e}chet differentiable with the Fr\'{e}chet derivative denoted by $F_i'(x)$ for
$x\in \D$. We will also assume that (\ref{sys}) has a solution. In general, (\ref{sys}) may
have many solutions. In order to find the desired one, some selection criteria should be enforced. We choose a proper,
lower semi-continuous, $p$-convex function $\Theta: \X\to (-\infty, \infty]$. By picking $x_0\in D(\partial \Theta)$
and $\xi_0\in \partial \Theta(x_0)$ as the initial guess, which may incorporate some available information on the sought
solution, we define $x^{\dag}$ to be the solution of (\ref{sys}) with the property
\begin{equation}\label{xdag}
\fl D_{\xi_0} \Theta(x^{\dag},x_0) := \min_{x\in D(\Theta)\cap \D }\left\{D_{\xi_0} \Theta(x,x_0) : F_i(x) = y_i,
i=0, \cdots, N-1\right\}.
\end{equation}

We will work under the following conditions on the operators $F_i$ where $B_\rho(x_0):=\{x\in \X: \|x-x_0\|\le \rho\}$.

\begin{Assumption}\label{A2}
\begin{enumerate}

\item[(a)] There is $\rho > 0$ such that $B_{2\rho}(x_0)\subset \D$ and (\ref{sys}) has a solution
in $B_\rho(x_0)\cap D(\Theta)$;

\item[(b)] Each operator $F_i$ is weakly closed on $\D$ and is Fr\'{e}chet differentiable on $B_{2\rho}(x_0)$,
and $x\to F_i'(x)$ is continuous
on $B_{2\rho}(x_0)$.

\item[(c)] Each $F_i$ is properly scaled so that $\|F_i'(x)\|\le 1$ for $x\in B_{2\rho}(x_0)$.

\item[(d)] There exists $0\le \eta<1$ such that
\begin{equation*}
\|F_i(x)-F_i(\bar x) -F_i'(\bar x) (x-\bar x)\|\le \eta \|F_i(x) -F_i(\bar x)\|
\end{equation*}
for all $x, \bar x \in B_{2\rho}(x_0)$ and $i=0, \cdots, N-1$.
\end{enumerate}
\end{Assumption}

All the conditions in Assumption \ref{A2} are standard. Condition (d) is called the tangential cone condition
and is widely used in the analysis of regularization methods for solving nonlinear ill-posed inverse problems (\cite{HNS95})
The weakly closedness of $F_i$ over $\D$ in (b) means that if $\{x_n\}\subset \D$ converges weakly to some $x\in \X$ and
$\{F_i(x_n)\}$ converges weakly to some $y_i\in \Y_i$, then $x\in \D$ and $F_i(x)=y_i$.

When $\X$ is a reflexive Banach space, by using the $p$-convexity and the weakly lower semi-continuity of $\Theta$ together with
the weakly closedness of $F_i$ for $i=0, \cdots, N-1$ it is standard to show that $x^\dag$ exists. The following result shows
that $x^\dag$ is in fact uniquely defined.

\begin{lemma}\label{lem0}
Let $\X$ be reflexive and $F_i$ satisfy Assumption \ref{A2}. If $x^\dag\in B_\rho(x_0)\cap D(\Theta)$, then $x^\dag$ is
the unique solution of (\ref{sys}) in $B_{2\rho}(x_0)\cap D(\Theta)$ satisfying (\ref{xdag}).
\end{lemma}

\proof{
Assume that (\ref{sys}) has another solution $\hat{x}$ in $B_{2\rho}(x_0)\cap D(\Theta)$ satisfying
(\ref{xdag}) with $\hat x\ne x^\dag$. Since $F_i(\hat{x})=y_i=F_i(x^\dag)$ for $i=0, \cdots, N-1$, we can use Assumption \ref{A2} (d)
to derive that
\begin{equation*}
F_i'(x^\dag) (\hat{x}-x^\dag)=0, \qquad i=0, \cdots, N-1.
\end{equation*}
Let $x_\la=\la \hat{x}+(1-\la) x^\dag$ for $0<\la<1$. Then $x_\la\in B_{2\rho}(x_0) \cap D(\Theta)$
and
\begin{equation*}
F_i'(x^\dag) (x_\la-x^\dag)=0, \qquad i=0, \cdots, N-1.
\end{equation*}
Thus we can use Assumption \ref{A2} (d) to conclude that
\begin{equation*}
\|F_i(x_\la)-F_i(x^\dag)\| \le \eta \|F_i(x_\la)-F_i(x^\dag)\|.
\end{equation*}
Since $0\le \eta <1$, this implies that $F_i(x_\la)=F_i(x^\dag)=y_i$ for $i=0, \cdots, N-1$.
Consequently, by the minimal property of $x^\dag$ we have
\begin{equation}\label{eq:20.7june}
D_{\xi_0} \Theta(x_\la, x_0) \ge  D_{\xi_0} \Theta(x^\dag, x_0).
\end{equation}
On the other hand, it follows from the strictly convexity of $\Theta$ that
\begin{eqnarray*}
\fl D_{\xi_0} \Theta (x_\la, x_0) %& = \Theta(x_\la)- \Theta(x_0) - \langle \xi_0, x_\la-x_0\rangle\\
%& < \la \Theta(\hat{x}) +(1-\la) \Theta(x^\dag) -\Theta(x_0) -\langle \xi_0, x_\la -x_0 \rangle \\
< \la D_{\xi_0} \Theta (\hat{x}, x_0) +(1-\la) D_{\xi_0} \Theta (x^\dag, x_0)
= D_{\xi_0} \Theta(x^{\dag}, x_0)
\end{eqnarray*}
for $0<\la<1$ which is a contradiction to (\ref{eq:20.7june}).
}

In practical application, instead of $y_i$ we only have noisy data $y_i^\d$ satisfying
\begin{equation*}
\|y_i^\d-y_i\| \le \d, \qquad i=0, \cdots, N-1
\end{equation*}
with a small known noise level $\d>0$. We will use $y_i^\d$, $i=0, \cdots, N-1$, to construct an approximate solution to (\ref{sys}).
We assume that each $\Y_i$ is uniformly smooth so that, for each $1<r<\infty$, the duality mapping $J_r^{\Y_i}: \Y_i\to \Y_i^*$
is single valued and continuous. By introducing a proper, lower-semi continuous, $p$-convex function $\Theta: (-\infty, \infty] \to \X$
satisfying (\ref{pconv}) for some constant $c_0>0$, we propose the following Landweber iteration of Kaczmarz type:

\begin{algorithm}\label{alg}

\begin{enumerate}

\item[(i)] Pick $\xi_0\in \X^*$ and set $x_0:=\arg\min_{x\in \X} \left\{\Theta(x)-\l \xi_0, x\r\right\}$;

\item[(ii)] Let $\xi_0^\d:=\xi_0$ and $x_0^\d:=x_0$. Assume that $\xi_n^\d$ and $x_n^\d$ are defined for some $n$,
we set $\xi_{n,0}^\d=\xi_n^\d$, $x_{n,0}^\d=x_n^\d$, and define
\begin{eqnarray*}
\xi_{n, i+1}^\d &= \xi_{n, i}^\d -\mu_{n, i}^\d F_i'(x_{n,i}^\d)^* J_r^{\Y_i}(F_i(x_{n,i}^\d)-y_i^\d),\\
x_{n, i+1}^\d &= \arg\min_{x\in \X} \left\{\Theta(x) -\l \xi_{n, i+1}^\d, x\r\right\}
\end{eqnarray*}
for $i=0, \cdots, N-1$, where
\begin{eqnarray*}
\mu_{n, i}^\d =\left\{\begin{array}{lll}
\mu_0 \|F_i(x_{n,i}^\d)-y_i^\d\|^{p-r}, &\quad  \mbox{ if } \|F_i(x_{n,i}^\d)-y_i^\d\| >\tau \d,\\
0, & \quad \mbox{ otherwise}
\end{array}\right.
\end{eqnarray*}
for some $\mu_0>0$. We then define $\xi_{n+1}^\d :=\xi_{n, N}^\d$ and $x_{n+1}^\d:=x_{n, N}^\d$.

\item[(iii)] Let $n_\d$ be the first integer such that
\begin{equation*}
\mu_{n_\d, i}^\d =0 \quad \mbox{for all } i=0, \cdots, N-1
\end{equation*}
and use $x_{n_\d}^\d$ to approximate the solution of (\ref{sys}).

\end{enumerate}

\end{algorithm}

In Algorithm \ref{alg}, $x_{n, i+1}^\d$ is defined as the minimizer of a $p$-convex functional over $\X$ which is independent of
$F_i$ and therefore it could be found by efficient solvers. By using (\ref{8.21.1}), one can see that
\begin{equation*}
x_{n, i+1}^\d=\nabla \Theta^*(\xi_{n, i+1}^\d)
\end{equation*}
which is useful for the forthcoming theoretical analysis.

In this section we will show that Algorithm \ref{alg} is well-defined by showing that $n_\d$ is finite and establish a
convergence result on $x_{n_\d}^\d$ as $\d\rightarrow 0$.

\begin{lemma}\label{L5.1}
Let Assumption \ref{A2} hold and let $\Theta: (-\infty, \infty] \to \X$ be a proper, lower-semi continuous, $p$-convex 
function with $p\ge 2$ satisfying (\ref{pconv}) for some $c_0>0$. Assume that
\begin{equation}\label{5.11.1}
D_{\xi_0} \Theta (x^\dag, x_0) \le c_0 \rho^p.
\end{equation}
Let $\{\xi_n^\d\}$ and $\{x_n^\d\}$ be defined by Algorithm \ref{alg} with
$\tau>1$ and $\mu_0>0$ such that
\begin{equation*}
c_1:=1-\eta-\frac{1+\eta}{\tau} - \frac{p-1}{p}\left(\frac{\mu_0}{2 c_0}\right)^{\frac{1}{p-1}}>0.
\end{equation*}
Then $n_\d<\infty$ and $x_{n,i}^\d\in B_{2\rho}(x_0)$ for all $n\ge 0$ and $i=0, \cdots, N-1$. Moreover,
for any solution $\hat x$ of (\ref{sys}) in $B_{2\rho}(x_0)\cap D(\Theta)$ and all $n$ there hold
\begin{eqnarray}
D_{\xi_{n+1}^\d}\Theta(\hat x, x_{n+1}^\d)  \le D_{\xi_n^\d}\Theta(\hat x, x_n^\d),\label{eq5.10}\\
\fl \,\,\,\, \quad \quad c_1\sum_{i=0}^{N-1} \mu_{n, i}^\d \|F_i(x_{n,i}^\d)-y_i^\d\|^r \le D_{\xi_n^\d}\Theta(\hat x, x_n^\d)
-D_{\xi_{n+1}^\d}\Theta(\hat x, x_{n+1}^\d). \label{eq5.11}
\end{eqnarray}
\end{lemma}

\proof{
In order to obtain (\ref{eq5.10}) and (\ref{eq5.11}), it suffices to show that $x_{n,i}^\d\in B_{2\rho}(x_0)$ and
\begin{equation}\label{12.16}
D_{\xi_{n, i+1}^\d}\Theta(\hat x, x_{n,i+1}^\d) -D_{\xi_{n,i}^\d}\Theta(\hat x, x_{n,i}^\d)
\le - c_1\mu_{n,i}^\d \|F_i(x_{n,i}^\d)-y_i^\d\|^r
\end{equation}
for all $n\ge 0$ and $i=0, \cdots, N-1$. From the definition of Bregman distance and (\ref{3.18.1}) it follows that
\begin{eqnarray*}
\fl D_{\xi_{n, i+1}^\d}\Theta(\hat x, x_{n,i+1}^\d) -D_{\xi_{n,i}^\d}\Theta(\hat x, x_{n,i}^\d)\\
= \Theta(x_{n, i}^\d) -\Theta(x_{n, i+1}^\d) -\l \xi_{n, i+1}^\d, \hat x -x_{n, i+1}^\d\r + \l \xi_{n, i}^\d, \hat x -x_{n, i}^\d\r\\
=\Theta^*(\xi_{n, i+1}^\d) -\Theta^*(\xi_{n, i}^\d) -\l \xi_{n, i+1}^\d- \xi_{n, i}^\d, \hat{x}\r.
\end{eqnarray*}
Using $x_{n,i}^\d=\nabla \Theta^*(\xi_{n,i}^\d)$, we can write
\begin{eqnarray*}
\fl D_{\xi_{n, i+1}^\d}\Theta(\hat x, x_{n,i+1}^\d) -D_{\xi_{n,i}^\d}\Theta(\hat x, x_{n,i}^\d)
&= \Theta^*(\xi_{n,i+1}^\d) -\Theta^*(\xi_{n,i}^\d) \\
&\quad \, -\l \xi_{n,i+1}^\d-\xi_{n,i}^\d, \nabla \Theta^*(\xi_{n,i}^\d)\r\\
& \quad \, +\l \xi_{n,i+1}^\d-\xi_{n,i}^\d, x_{n,i}^\d-\hat{x}\r.
\end{eqnarray*}
Since $\Theta$ is $p$-convex, we may use (\ref{6.21.1}) to obtain
\begin{eqnarray*}
\fl D_{\xi_{n, i+1}^\d}\Theta(\hat x, x_{n,i+1}^\d) -D_{\xi_{n,i}^\d}\Theta(\hat x, x_{n,i}^\d)
 \le \frac{1}{p^*(2 c_0)^{p^*-1}} \|\xi_{n, i+1}^\d-\xi_{n,i}^\d\|^{p^*} \\
 -\mu_{n,i}^\d \l J_r^{\Y_i} (F_i(x_{n,i}^\d)-y_i^\d), F_i'(x_{n,i}^\d)(x_{n,i}^\d-\hat{x})\r,
\end{eqnarray*}
where $1/p+1/p^*=1$. By using the properties of the duality mapping $J_r^{\Y_i}$ and Assumption \ref{A2} it follows that
\begin{eqnarray}\label{12.16.1}
\fl D_{\xi_{n, i+1}^\d} \Theta(\hat x, x_{n,i+1}^\d) -D_{\xi_{n,i}^\d}\Theta(\hat x, x_{n,i}^\d)\nonumber\\
 \le  \mu_{n,i}^\d \| F_i(x_{n,i}^\d)-y_i^\d\|^{r-1}\|y_i^\d-F_i(x_{n,i}^\d)-F_i'(x_{n,i}^\d)(\hat{x}-x_{n,i}^\d)\| \nonumber\\
  \quad -\mu_{n,i}^\d \|F_i(x_{n,i}^\d)-y_i^\d\|^r+ \frac{1}{p^*(2 c_0)^{p^*-1}}\|\xi_{n, i+1}^\d-\xi_{n,i}^\d\|^{p^*} \nonumber\\
\le (1+\eta) \mu_{n,i}^\d \|F_i(x_{n,i}^\d)-y_i^\d\|^{r-1} \d - (1-\eta) \mu_{n,i}^\d \|F_i(x_{n,i}^\d)-y_i^\d\|^r \nonumber\\
  \quad + \frac{1}{p^*(2 c_0)^{p^*-1}}\left(\mu_{n,i}^\d\right)^{p^*} \|F_i'(x_{n,i}^\d)^* J_r^{\Y_i} (F_i(x_{n,i}^\d)-y_i^\d)\|^{p^*}.
\end{eqnarray}
According to the definition of $\mu_{n,i}^\d$, the scaling condition in Assumption \ref{A2} (c), and the property of $J_r^{\Y_i}$, it is easy to see that
\begin{eqnarray*}
\quad \mu_{n,i}^\d \|F_i(x_{n,i}^\d)-y_i^\d\|^{r-1} \d &\le \frac{1}{\tau} \mu_{n,i}^\d \|F_i(x_{n,i}^\d)-y_i^\d\|^r,\\
\fl \,\,\left(\mu_{n,i}^\d\right)^{p^*-1} \|F_i'(x_{n,i}^\d)^* J_r^{\Y_i} (F_i(x_{n,i}^\d)-y_i^\d)\|^{p^*} &\le \mu_0^{p^*-1} \|F_i(x_{n,i}^\d)-y_i^\d\|^r.
\end{eqnarray*}
Combining these two inequalities with (\ref{12.16.1}) we can obtain (\ref{12.16}). To show $x_{n,i+1}^\d \in B_{2\rho}(x_0)$
we first use (\ref{12.16}) with $\hat x =x^\dag$ and (\ref{5.11.1}) to obtain
\begin{equation*}
D_{\xi_{n,i+1}^\d}\Theta(x^\dag, x_{n, i+1}^\d) \le D_{\xi_0}\Theta(x^\dag, x_0) \le c_0 \rho^p.
\end{equation*}
In view of (\ref{pconv}), we then have $\|x_{n, i+1}^\d-x^\dag\|\le \rho$ and $\|x^\dag-x_0\|\le \rho$.
Consequently $x_{n, i+1}^\d\in B_{2\rho}(x_0)$.

We next show $n_\d<\infty$. According to the definition of $n_\d$, for any $n<n_\d$ there is at least one
$i_n\in \{0, \cdots, N-1\}$ such that $\|F_{i_n}(x_{n, i_n}^\d)-y_{i_n}^\d\|>\tau \d$. Consequently
\begin{equation*}
\mu_{n, i_n}^\d =\mu_0 \|F_{i_n}(x_{n, i_n}^\d)-y_{i_n}^\d\|^{p-r}
\end{equation*}
and
\begin{equation*}
\sum_{i=0}^{N-1} \mu_{n,i}^\d \|F_i(x_{n,i}^\d)-y_i^\d\|^r \ge \mu_0 \|F_{i_n}(x_{n, i_n}^\d)-y_{i_n}^\d\|^p
>\mu_0 \tau^p \d^p.
\end{equation*}
By summing (\ref{eq5.11}) over $n$ from $n=0$ to $n=m$ for any $m<n_\d$ and using the above inequality we obtain
$c_1 \mu_0 \tau^p \d^p (m+1) \le D_{\xi_0}\Theta(\hat x, x_0)$. Since this is true for any $m<n_\d$, we must have $n_\d<\infty$.
}

When using the exact data $y_i$ instead of the noisy data $y_i^\d$ in Algorithm \ref{alg}, we will drop
the superscript $\d$ in all the quantities involved, for instance, we will write $\xi_n^\d$ as $\xi_n$,
$x_n^\d$ as $x_n$, and so on. Observing that
\begin{equation*}
\mu_{n,i} \|F_i(x_{n,i})-y_i\|^r =\mu_0 \|F_i(x_{n,i})-y_i\|^p.
\end{equation*}
The proof of Lemma \ref{L5.1} in fact shows that, under Assumption \ref{A2}, if
\begin{equation*}
c_2:=1-\eta -\frac{p-1}{p}\left(\frac{\mu_0}{2 c_0}\right)^{\frac{1}{p-1}}>0,
\end{equation*}
then
\begin{equation*}
x_{n,i}\in B_{2\rho}(x_0) \qquad \forall n\ge 0 \mbox{ and } i=0, \cdots, N-1
\end{equation*}
and for any solution $\hat x$ of (\ref{sys}) in $B_{2\rho}(x_0)\cap D(\Theta)$ and all $n$ there hold
\begin{eqnarray}
D_{\xi_{n+1}}\Theta(\hat x, x_{n+1})  \le D_{\xi_n}\Theta(\hat x, x_n),\label{eq5.10.1}\\
c_2\mu_0 \sum_{i=0}^{N-1} \|F_i(x_{n,i})-y_i\|^p  \le D_{\xi_n}\Theta(\hat x, x_n)
-D_{\xi_{n+1}}\Theta(\hat x, x_{n+1}). \label{eq5.11.1}
\end{eqnarray}
These two inequalities imply immediately that
\begin{equation}\label{eq5.12.1}
\lim_{n\rightarrow \infty} \sum_{i=0}^{N-1} \|F_i(x_{n,i})-y_i\|^p =0.
\end{equation}
The next result gives an estimate on $\|F_i(x_n)-y_i\|$ and shows that $\|F_i(x_n)-y_i\|\rightarrow 0$ as
$n\rightarrow \infty$ for all $i=0, \cdots, N-1$.

\begin{lemma}\label{L5.9}
Let all the conditions in Lemma \ref{L5.1} hold. Then there is a constant $C_0$ such that for all $n\ge 1$ there hold
\begin{equation*}
\|F_i(x_n)-y_i\| \le C_0 \sum_{j=0}^i \|F_j(x_{n,j})-y_j\|, \quad i=0, \cdots, N-1.
\end{equation*}
Consequently $\lim_{n\rightarrow \infty} \|F_i(x_n)-y_i\|=0$ for all $i=0, \cdots, N-1$.
\end{lemma}

\proof{
Recall that $x_n=x_{n, 0}$, we have
\begin{eqnarray*}
\|F_i(x_n)-y_i\|\le \|F_i(x_{n, i})-y_i\| +\|F_i(x_{n,i})-F_i(x_{n, 0})\|.
\end{eqnarray*}
By using the condition on $F$ we have
\begin{eqnarray}\label{12.15.1}
\fl \quad \,\, \|F_i(x_n)-y_i\| &\le \|F_i(x_{n, i})-y_i\| + \sum_{j=0}^{i-1} \|F_i(x_{n,j+1})-F_i(x_{n,j})\| \nonumber\\
& \le  \|F_i(x_{n, i})-y_i\| +\frac{1}{1-\eta} \sum_{j=0}^{i-1} \|F_i'(x_{n,j})(x_{n,j+1}-x_{n,j})\| \nonumber\\
& \le  \|F_i(x_{n, i})-y_i\| + \frac{1}{1-\eta} \sum_{j=0}^{i-1} \|x_{n,j+1}-x_{n,j}\|.
\end{eqnarray}
Since $x_{n,j}=\nabla \Theta^*(\xi_{n,j})$, we can use the property (\ref{3.29.1}) to derive that
\begin{eqnarray*}
 \|x_{n,j+1}-x_{n,j}\| &=\|\nabla \Theta^*(\xi_{n,j+1})-\nabla \Theta^*(\xi_{n,j})\|\\
& \le \left(\frac{\|\xi_{n,j+1} -\xi_{n,j}\|}{2 c_0}\right)^{\frac{1}{p-1}}.
\end{eqnarray*}
Using the definition of $\xi_{n,j+1}$ and the property of the duality mapping $J_r^{\Y_i}$ we obtain
\begin{eqnarray*}
\|x_{n,j+1}-x_{n,j}\| &\le (2c_0)^{-\frac{1}{p-1}} \mu_{n,j}^{\frac{1}{p-1}} \|F_j(x_{n,j}) -y_j\|^{\frac{r-1}{p-1}} \\
&=\left(\frac{\mu_0}{2 c_0}\right)^{\frac{1}{p-1}}  \|F_j(x_{n,j})-y_j\|.
\end{eqnarray*}
Combining this with (\ref{12.15.1}) gives the desired inequality.
}

As the first step toward the proof of convergence on $x_{n_\d}^\d$, we need to derive some convergence results on the
sequences $\{x_n\}$ and $\{\xi_n\}$.  This will be achieved by the following proposition which gives a general convergence
criterion on any sequences $\{x_n\}\subset \X$ and $\{\xi_n\}\subset \X^*$ satisfying certain conditions.

\begin{proposition}\label{general}
Consider the system (\ref{sys}) for which Assumption \ref{A2} holds.
Let $\Theta: \X\to (-\infty, \infty]$ be a proper, lower semi-continuous and uniformly convex function.
Let $\{x_n\}\subset B_{2\rho}(x_0)$ and $\{\xi_n\}\subset \X^*$ be such that

\begin{enumerate}
\item[(i)] $\xi_n\in \p \Theta(x_n)$ for all $n$;

\item[(ii)] for any solution $\hat x$ of (\ref{sys}) in $B_{2\rho}(x_0)\cap D(\Theta)$ the sequence
$\{D_{\xi_n}\Theta(\hat x, x_n)\}$ is monotonically decreasing;

\item[(iii)] $\lim_{n\rightarrow \infty} \|F_i(x_n)-y_i\|=0$ for all $i=0, \cdots, N-1$.

\item[(iv)] there is a subsequence $\{n_k\}$ with $n_k\rightarrow \infty$ such that for all $l<k$ and any solution
$\hat x$ of (\ref{sys}) in $B_{2\rho}(x_0)\cap D(\Theta)$ there holds
\begin{equation}\label{12.15.11}
\fl \qquad \qquad |\l \xi_{n_k}-\xi_{n_l}, x_{n_k}-\hat x\r|\le C_1\left(D_{\xi_{n_l}}\Theta(\hat x, x_{n_l}) -D_{\xi_{n_k}}\Theta(\hat x, x_{n_k})\right)
\end{equation}
for some constant $C_1$.
\end{enumerate}
Then there exists a solution $x_*$ of (\ref{sys}) in $B_{2\rho}(x_0)\cap D(\Theta)$ such that
\begin{equation*}
\lim_{n\rightarrow \infty} D_{\xi_n}\Theta(x_*, x_n)=0.
\end{equation*}
If, in addition, $x^\dag \in B_\rho(x_0)\cap D(\Theta)$ and $\xi_{n+1}-\xi_n \in \overline{{\mathcal R}(F_0'(x^\dag)^*)}+\cdots + \overline{{\mathcal R}(F_{N-1}'(x^\dag)^*)}$
for all $n$, then $x_*=x^\dag$.
\end{proposition}

Proposition \ref{general} may be of independent interest and its proof is deferred to section \ref{P3}.
Now we are ready to give the convergence result on $\{\xi_n\}$ and $\{x_n\}$ defined by the Landweber
iteration of Kaczmarz type with exact data.

\begin{lemma}\label{L5.21}
Let all the conditions in Lemma \ref{L5.1} hold. For the sequences $\{\xi_n\}$ and $\{x_n\}$ defined by Algorithm \ref{alg} with
exact data, there exists a solution $x_*\in B_{2\rho}(x_0)\cap D(\Theta)$ of (\ref{sys}) such that
\begin{equation*}
\lim_{n\rightarrow \infty} \|x_n-x_*\|=0 \qquad \mbox{and} \qquad \lim_{n\rightarrow \infty} D_{\xi_n}\Theta(x_*, x_n)=0.
\end{equation*}
If in addition $\N(F_i'(x^\dag))\subset \N(F_i'(x))$ for all $x\in B_{2\rho}(x_0)$ and $i=0, \cdots, N-1$, then $x_*=x^\dag$.
\end{lemma}

\proof{
We will use Proposition \ref{general} to complete the proof. By the definition of $\{\xi_n\}$ and $\{x_n\}$ we have
$\xi_n\in \p \Theta(x_n)$. The monotonicity of $\{D_{\xi_n}\Theta(\hat x, x_n)\}$ is given by (\ref{eq5.10.1}). Lemma \ref{L5.9}
shows that
\begin{equation*}
\lim_{n\rightarrow \infty} \|F_i(x_n)-y_i\|=0 \qquad \mbox{ for all } i=0, \cdots, N-1.
\end{equation*}
Therefore, in order to derive the convergence result, it suffices to show that there exists a strictly increasing subsequence
$\{n_k\}$ such that for any solution $\hat{x}$ of (\ref{sys}) and any $l<k$ there holds
\begin{equation}\label{eq5.4}
\left|\l \xi_{n_k}-\xi_{n_l}, x_{n_k}-\hat x\r\right| \le C\left(D_{\xi_{n_l}}\Theta(\hat x, x_{n_l}) -D_{\xi_{n_k}}\Theta (\hat x, x_{n_k})\right)
\end{equation}
To this end, let
\begin{equation*}
R_n:=\sum_{i=0}^{N-1} \|y_i-F_i(x_{n, i})\|^p.
\end{equation*}
It follows from (\ref{eq5.12.1}) that
\begin{equation}\label{eq5.1}
\lim_{n\rightarrow \infty} R_n =0.
\end{equation}
Moreover, if $R_n=0$ for some $n$, then $y_i=F_i(x_{n,i})$ for $i=0, \cdots, N-1$. Consequently it follows from the
definition of the method that $x_{m,i}=x_n$ for all $m\ge n$ and $i=0, \cdots, N-1$. Therefore
\begin{equation}\label{eq5.2}
R_n=0 \mbox{ for some } n \Longrightarrow R_m=0 \mbox{ for all } m\ge n.
\end{equation}
In view of (\ref{eq5.1}) and (\ref{eq5.2}), we can introduce a subsequence $\{n_k\}$ by setting $n_0=0$ and letting $n_k$,
for each $k\ge 1$, be the first integer satisfying
\begin{equation*}
n_k\ge n_{k-1}+1 \qquad \mbox{and} \qquad R_{n_k} \le R_{n_{k-1}}.
\end{equation*}
For such chosen strictly increasing sequence $\{n_k\}$ it is easy to see that
\begin{equation}\label{eq5.3}
R_{n_k} \le R_n, \qquad 0\le n<n_k.
\end{equation}

We now prove (\ref{eq5.4}) for the above chosen subsequence $\{n_k\}$. We first use the definition of the method
to obtain for $n<n_k$ that
\begin{eqnarray*}
\fl \l \xi_{n+1}-\xi_n, x_{n_k}-\hat x\r & =\sum_{i=0}^{N-1} \l \xi_{n, i+1}-\xi_{n,i}, x_{n_k}-\hat x\r\\
&=-\sum_{i=0}^{N-1} \mu_{n, i} \l J_r^{\Y_i} (F_i(x_{n,i})-y_i), F_i'(x_{n,i})(x_{n_k}-\hat x)\r.
\end{eqnarray*}
Using the condition on $F_i$ it is easy to obtain
\begin{eqnarray*}
\fl \|F_i'(x_{n,i})(x_{n_k}-\hat x)\| &\le \|F_i'(x_{n,i})(x_{n,i}-\hat x)\| +\|F_i'(x_{n,i})(x_{n_k}-x_{n, i})\|\\
& \le (1+\eta) \left(\|F_i(x_{n,i})-y_i\|+ \|F_i(x_{n_k})-F_i(x_{n,i})\|\right)\\
&\le (1+\eta) \left(2 \|F_i(x_{n,i})-y_i\|+\|F_i(x_{n_k})-y_i\|\right).
\end{eqnarray*}
Therefore, by using the property of the duality mapping $J_r^{\Y_i}$ and the H\"{o}lder inequality, we have
\begin{eqnarray}\label{eq5.6}
\fl |\l \xi_{n+1}-\xi_n, x_{n_k}-\hat x\r| \le \sum_{i=0}^{N-1} \mu_{n, i} \|F_i(x_{n,i})-y_i\|^{r-1}
\|F_i'(x_{n,i})(x_{n_k}-\hat x)\| \nonumber\\
\le (1+\eta) \sum_{i=0}^{N-1} \mu_{n,i} \|F_i(x_{n,i})-y_i\|^{r-1} \|F_i(x_{n_k})-y_i\|+2(1+\eta) \mu_0 R_n \nonumber\\
= (1+\eta) \mu_0 \sum_{i=0}^{N-1} \|F_i(x_{n,i})-y_i\|^{p-1} \|F_i(x_{n_k})-y_i\|+2(1+\eta) \mu_0 R_n \nonumber\\
\le (1+\eta) \mu_0 R_n^{\frac{p-1}{p}} \left(\sum_{i=0}^{N-1} \|F_i(x_{n_k})-y_i\|^p\right)^{\frac{1}{p}} +2(1+\eta)\mu_0 R_n.
\end{eqnarray}
From Lemma \ref{L5.9} and the H\"{o}lder inequality it follows that
\begin{eqnarray*}
\fl \|F_i(x_{n_k})-y_i\| & \le  C_0 \sum_{j=0}^i \|F_j(x_{n_k,j})-y_j\|
 \le  C_0 (i+1)^{\frac{p-1}{p}} \left(\sum_{j=0}^i \|F_j(x_{n_k,j})-y_j\|^p\right)^{\frac{1}{p}}.
\end{eqnarray*}
This implies that
\begin{eqnarray*}
\sum_{i=0}^{N-1} \|F_i(x_{n_k})-y_i\|^p &\le C_0^p \sum_{i=0}^{N-1} \sum_{j=0}^i (i+1)^{p-1} \|F_j(x_{n_k,j})-y_j\|^p\\
& \le C_0^p (N+1)^p \sum_{j=0}^{N-1} \|F_j(x_{n_k,j}) -y_j\|^p \\
& =C_0^p (N+1)^p R_{n_k}.
\end{eqnarray*}
Combining this with (\ref{eq5.6}) and using (\ref{eq5.3}) we can obtain for $n<n_k$ that
\begin{eqnarray*}
|\l \xi_{n+1}-\xi_n, x_{n_k}-\hat x\r| &\le (1+\eta)\mu_0 \left[C_0(N+1) R_n^{\frac{p-1}{p}} R_{n_k}^{\frac{1}{p}} + 2 R_n\right] \\
& \le (1+\eta) \mu_0 \left[C_0(N+1)+2\right] R_n.
\end{eqnarray*}
Finally, we can derive that
\begin{eqnarray*}
\left|\l \xi_{n_k}-\xi_{n_l}, x_{n_k}-\hat x\r\right|\le\sum_{n=n_l}^{n_k-1} |\l \xi_{n+1}-\xi_n, x_{n_k}-\hat x\r|
\le C_2 \sum_{n=n_l}^{n_k-1} R_n
\end{eqnarray*}
with $C_2:=(1+\eta) \mu_0 \left[C_0(N+1)+2\right]$. In view of (\ref{eq5.11.1}) we therefore obtain (\ref{eq5.4}).

In order to show $x_*=x^\dag$ under the condition $\N(F_i'(x^\dag))\subset \N(F_i'(x))$ for all
$x\in B_{2\rho}(x_0)$ and $i=0, \cdots, N-1$, we observe from the definition of $\xi_n$ that
\begin{equation*}
\fl \xi_{n+1}-\xi_n =\sum_{i=0}^{N-1} (\xi_{n, i+1} -\xi_{n, i})\in {\mathcal R}(F_0'(x_{n,0})^*)+\cdots + {\mathcal R}(F_{N-1}'(x_{n, N-1})^*).
\end{equation*}
Since $\X$ is reflexive and $\N(F_i'(x^\dag)) \subset \N(F_i'(x_{n,i}))$, we have
$\overline{\R(F_i'(x_{n,i})^*)}\subset \overline{\R(F_i'(x^\dag)^*)}$. Therefore
\begin{equation*}
\xi_{n+1}-\xi_n \in \overline{{\mathcal R}(F_0'(x^\dag)^*)}+\cdots + \overline{{\mathcal R}(F_{N-1}'(x^\dag)^*)}
\quad \mbox{ for all } n.
\end{equation*}
Hence, we can use the second part of Proposition \ref{general} to conclude that $x_*=x^\dag$.
}

In order to use the above result to prove the convergence of the Landweber iteration of Kaczmarz type
described in Algorithm \ref{alg} with noisy data, we need the following stability result.

\begin{lemma}\label{stability}
Let $\X$ be reflexive and let $\Y_i$ be uniformly smooth. Let all the conditions in Lemma \ref{L5.1} hold. Then for all $n\ge 0$ and $i=0, \cdots, N-1$ there hold
\begin{equation*}
\xi_{n, i}^\d \rightarrow \xi_{n, i} \quad \mbox{ and } \quad x_{n,i}^\d\rightarrow x_{n,i} \quad \mbox{as } \d\rightarrow 0.
\end{equation*}
\end{lemma}

\proof{
The result is trivial for $n=0$ and $i=0$. We next assume that the result is true for some $n\ge 0$ and
some $i\in \{0, \cdots, N-1\}$ and show that $\xi_{n, i+1}^\d\rightarrow \xi_{n, i+1}$ and $x_{n, i+1}^\d\rightarrow x_{n,i+1}$
as $\d\rightarrow 0$. We consider two cases.

{\it Case 1: $F_i(x_{n,i})=y_i$}. In this case we have $\mu_{n,i}=0$ and $\|F_i(x_{n,i}^\d) -y_i^\d\|\rightarrow 0$
as $\d\rightarrow 0$ by the continuity of $F_i$. Thus
\begin{equation*}
\xi_{n, i+1}^\d-\xi_{n,i+1} =\xi_{n,i}^\d-\xi_{n,i} -\mu_{n,i}^\d F_i'(x_{n,i}^\d)^* J_r^{\Y_i} (F_i(x_{n,i}^\d)-y_i^\d)
\end{equation*}
which implies that
\begin{equation*}
\|\xi_{n, i+1}^\d-\xi_{n,i+1}\| \le\|\xi_{n,i}^\d-\xi_{n,i}\| +\mu_0 \|F_i(x_{n,i}^\d)-y_i^\d\|^{p-1}.
\end{equation*}
By the induction hypotheses, we then have $\xi_{n,i+1}^\d\rightarrow \xi_{n, i+1}$ as $\d\rightarrow 0$. Consequently,
by using the continuity of $\nabla \Theta^*$ we have $x_{n, i+1}^\d =\nabla \Theta^*(\xi_{n, i+1}^\d)\rightarrow
\nabla \Theta^*(\xi_{n, i+1}) =x_{n, i+1}$ as $\d\rightarrow 0$.

{\it Case 2: $F_i(x_{n, i})\ne y_i$}. In this case we have $\|F_i(x_{n, i}^\d)-y_i^\d\|>\tau \d$ for small $\d>0$. Therefore
\begin{equation*}
\mu_{n, i}^\d=\mu_0 \|F_i(x_{n,i}^\d)-y_i^\d\|^{p-r} \rightarrow \mu_{n,i}=\mu_0 \|F_i(x_{n,i})-y_i\|^{p-r}
\end{equation*}
as $\d\rightarrow 0$. By Assumption \ref{A2} (b) and the uniform smoothness of $\Y_i$, we know that $F$, $F'$ and $J_r^{\Y_i}$
are continuous. It then follows from the induction hypotheses that $\xi_{n, i+1}^\d\rightarrow \xi_{n, i+1}$ and hence
$x_{n, i+1}^\d\rightarrow x_{n, i+1}$ as $\d\rightarrow 0$ using again the continuity of $\nabla \Theta^*$.
}

We are now in a position to give the main convergence result on the Landweber iteration of Kaczmarz type.

\begin{theorem}\label{T6.4}
Let $\X$ be reflexive and let $\Y_i$ be uniformly smooth, let Assumption \ref{A2} hold with
$0\le \eta<1$, and let $\Theta: \X\to (-\infty, \infty]$ be proper, lower semi-continuous, and
$p$-convex function satisfies (\ref{pconv}). Assume that (\ref{5.11.1}) holds. Then for $\{\xi_n^\d\}$ and $\{x_n^\d\}$
defined by Algorithm \ref{alg} with $\tau>1$ and $\mu_0>0$ satisfying
\begin{equation}\label{6.21.6}
1-\eta-\frac{1+\eta}{\tau}-\frac{p-1}{p}\left(\frac{\mu_0}{2 c_0}\right)^{\frac{1}{p-1}}>0,
\end{equation}
there is a solution $x_*\in B_{2\rho}(x_0)\cap D(\Theta)$ of (\ref{sys}) such that
\begin{equation*}
\lim_{\d\rightarrow 0} \|x_{n_\d}^\d-x_*\|=0 \qquad \mbox{and} \qquad
\lim_{\d\rightarrow 0} D_{\xi_{n_\d}^\d}\Theta(x_*, x_{n_\d}^\d) =0
\end{equation*}
If in addition $\N(F_i'(x^\dag))\subset \N(F_i'(x))$ for all $x\in B_{2\rho}(x_0)\cap D(F)$ and $i=0, \cdots, N-1$,
then $x_*=x^\dag$.
\end{theorem}
\proof{
Let $x_*$ be the solution of (\ref{sys}) determined in Lemma \ref{L5.21}. Due to the $p$-convexity of $\Theta$,
it suffices to show that $\lim_{\d\rightarrow 0} D_{\xi_{n_\d}^\d} \Theta(x_*, x_{n_\d}^\d)=0$.
We complete the proof by considering two cases.

Assume first that $\{y_i^{\d_k}\}$, $i=0, \cdots, N-1$, are sequences satisfying $\|y_i^{\d_k}-y_i\|\le \d_k$
with $\d_k\rightarrow 0$ such that $n_k:=n_{\d_k}\rightarrow \hat n$ as $k\rightarrow \infty$ for some finite
integer $\hat n$. We may assume $n_k=\hat n$ for all $k$. From the definition of $\hat n:=n_k$ we have
\begin{equation}\label{eq101}
\|F_i(x_{\hat n}^{\d_k})-y_i^{\d_k}\|\le \tau \d_k, \qquad i=0, \cdots, N-1.
\end{equation}
By taking $k\rightarrow \infty$ and using Lemma \ref{stability}, we can obtain
\begin{equation*}
F_i(x_{\hat n})=y_i, \qquad i=0, \cdots, N-1.
\end{equation*}
Using the definition of $\{\xi_n\}$ and $\{x_n\}$, this implies that $\xi_n=\xi_{\hat n}$ and
$x_n=x_{\hat n}$ for all $n\ge \hat n$. Since Lemma \ref{L5.21} implies that
$x_n\rightarrow x_*$ as $n\rightarrow \infty$, we must have $x_{\hat n}=x_*$. Consequently, by Lemma \ref{stability},
$\xi_{n_k}^{\d_k} \rightarrow \xi_{\hat n}$ and $x_{n_k}^{\d_k} \rightarrow x_*$ as $k\rightarrow \infty$.
This together with the lower semi-continuity of $\Theta$ implies that
\begin{eqnarray*}
0 & \le \liminf_{k\rightarrow \infty} D_{\xi_{n_k}^{\d_k}}\Theta(x_*, x_{n_k}^{\d_k})
\le \limsup_{k\rightarrow \infty} D_{\xi_{n_k}^{\d_k}}\Theta(x_*, x_{n_k}^{\d_k})\\
& = \Theta(x_*) -\liminf_{k\rightarrow \infty} \Theta(x_{n_k}^{\d_k})
- \lim_{k\rightarrow \infty} \l \xi_{n_k}^{\d_k}, x_*-x_{n_k}^{\d_k}\r\\
& \le \Theta(x_*)-\Theta(x_*)=0.
\end{eqnarray*}
Therefore  $\lim_{k\rightarrow \infty} D_{\xi_{n_k}^{\d_k}}\Theta(x_*, x_{n_k}^{\d_k})=0$.

Assume next that $\{y_i^{\d_k}\}$, $i=0, \cdots, N-1$, are sequences satisfying $\|y_i^{\d_k}-y_i\|\le \d_k$
with $\d_k\rightarrow 0$ such that $n_k:=n_{\d_k}\rightarrow \infty$ as $k\rightarrow \infty$.
Let $n$ be any fixed integer. then $n_k>n$ for large $k$. It then follows from (\ref{eq5.10}) in Lemma \ref{L5.1} that
\begin{eqnarray*}
D_{\xi_{n_k}^{\d_k}} \Theta(x_*, x_{n_k}^{\d_k}) \le D_{\xi_n^{\d_k}} \Theta (x_*, x_n^{\d_k})
=\Theta(x_*) -\Theta(x_n^{\d_k}) -\l \xi_n^{\d_k}, x_*-x_n^{\d_k}\r.
\end{eqnarray*}
By using Lemma \ref{stability} and the lower semi-continuity of $\Theta$ we obtain
\begin{eqnarray*}
0 & \le \liminf_{k\rightarrow \infty} D_{\xi_{n_k}^{\d_k}} \Theta (x_*, x_{n_k}^{\d_k})
\le \limsup_{k\rightarrow \infty} D_{\xi_{n_k}^{\d_k}} \Theta (x_*, x_{n_k}^{\d_k})\\
& \le \Theta(x_*) -\liminf_{k\rightarrow \infty} \Theta(x_n^{\d_k}) -\lim_{k\rightarrow \infty} \l \xi_n^{\d_k}, x_*-x_n^{\d_k}\r \\
& \le \Theta(x_*) -\Theta(x_n) -\l \xi_n, x_*-x_n\r\\
&= D_{\xi_n}\Theta(x_*, x_n).
\end{eqnarray*}
Since $n$ can be arbitrary and since Lemma \ref{L5.21} implies that $D_{\xi_n}\Theta(x_*, x_n)\rightarrow 0$
as $n\rightarrow \infty$, we therefore have $\lim_{k\rightarrow \infty} D_{\xi_{n_k}^{\d_k}} \Theta (x_*, x_{n_k}^{\d_k})=0$.
}

\begin{remark}\label{rem:6.21.1}
{\rm
Without assuming the scaling condition (c) in Assumption \ref{A2}, a good choice of the step length $\mu_{n,i}^\d$ in
Algorithm \ref{alg} requires the knowledge of $\|F_i'(x_{n,i}^\d)\|$ which is not easy to estimate. In order to
avoid this inconvenience, we may consider the alternative choice
\begin{eqnarray*}
\mu_{n, i}^\d =\left\{\begin{array}{lll}
\tilde{\mu}_{n, i}^\d \|F_i(x_{n, i}^\d)-y_i^\d\|^{p-r} &\quad  \mbox{ if } \|F_i(x_{n,i}^\d)-y_i^\d\| >\tau \d,\\
0, & \quad \mbox{ otherwise},
\end{array}\right.
\end{eqnarray*}
where
\begin{equation*}
\tilde{\mu}_{n, i}^\d: =\min\left\{\frac{\mu_0\|F_i(x_{n,i}^\d)-y_i^\d\|^{p(r-1)}}{\|F_i'(x_{n, i}^\d)^* J_r^{\Y_i} (F_i(x_{n, i}^\d)-y_i^\d)\|^p},
 \, \mu_1 \right\}
\end{equation*}
with two positive constants $\mu_0$ and $\mu_1$. It is easy to see that
\begin{equation*}
\min\{\mu_0 B_0^{-p}, \mu_1\} \le \tilde{\mu}_{n, i}^\d \le \mu_1,
\end{equation*}
where $B_0>0$ is a constant such that $\|F_i'(x)\|\le B_0$ for all $x\in B_{2\rho}(x_0)$ and $i=0, \cdots, N-1$.
If $\tau>1$ and $\mu_0>0$ are chosen such that (\ref{6.21.6}) holds and $\mu_1>0$ is chosen to be any number,
then, with some obvious modification in the proof of Lemma \ref{L5.1}, we can obtain
$n_\d<\infty$ and the monotonicity result. The requirement $\tilde{\mu}_{n,i}^\d\le \mu_1$ is to guarantee that the stability result in
Lemma \ref{stability} remains true. Therefore, we can still obtain the same convergence result as in Theorem \ref{T6.4}.
In order to allow large step lengths, we usually take $\mu_1$ to be a large number in practical applications.
}
\end{remark}

\begin{remark}
{\rm In Algorithm \ref{alg} we may replace the step (iii) by the stopping criterion that defines $n_\d$ to be the first integer
satisfying
\begin{equation*}
\sum_{i=0}^{N-1} \|F_i(x_{n_\d, i}^\d)-y_i^\d\|^p \le N \tau^p \d^p.
\end{equation*}
With some minor changes in the above arguments, we can still have the same convergence result as in Theorem \ref{T6.4}.
This new stopping criterion clearly terminates the iteration earlier than (iii) of Algorithm \ref{alg} and hence provides
an opportunity of avoiding computing many additional iterations that do not have essential contribution in the final stage. }
\end{remark}

\section{Proof of Proposition \ref{general}}\label{P3}

Proposition \ref{general} plays an important role in the convergence analysis of Algorithm \ref{alg} in Section 3.
It shows that for some sequences $\{x_n\}\subset \X$ and $\{\xi_n\}\subset \X^*$ constructed by a suitable algorithm,
the convergence of $\{x_n\}$ can be derived by showing certain monotonicity result together with a result like
(\ref{12.15.11}) along a suitable chosen subsequence of integers. This result might be useful for analyzing other methods
as well. In the following we give the proof.

\vskip 0.2cm

\noindent{\bf Proof of Proposition \ref{general}.}
We first show the convergence of $\{x_{n_k}\}$. For any $l<k$ we have from (\ref{4.3.1}) and (\ref{12.15.11}) that
\begin{eqnarray*}
\fl D_{\xi_{n_l}}\Theta(x_{n_k}, x_{n_l}) & =D_{\xi_{n_l}}\Theta(\hat{x}, x_{n_l}) -D_{\xi_{n_k}}\Theta(\hat{x}, x_{n_k})
+\l \xi_{n_k}-\xi_{n_l}, x_{n_k}-\hat{x}\r\\
&\le (1+C_1) \left(D_{\xi_{n_l}}\Theta(\hat{x}, x_{n_l}) -D_{\xi_{n_k}}\Theta(\hat{x}, x_{n_k})\right).
\end{eqnarray*}
By the monotonicity of $\{D_{\xi_n}\Theta(\hat x, x_n)\}$ we can conclude that $D_{\xi_{n_l}}\Theta(x_{n_k}, x_{n_l})\rightarrow 0$
as $k, l\rightarrow \infty$. In view of the uniformly convexity of $\Theta$, it follows that $\{x_{n_k}\}$ is a
Cauchy sequence in $\X$. Thus $x_{n_k}\rightarrow x_*$ for some $x_*\in B_{2\rho}(x_0)\subset \X$. Since $\lim_{n\rightarrow \infty} \|F_i(x_n) -y_i\|=0$,
we have $F_i(x_*)=y_i$ for all $i=0, \cdots, N-1$.

In order to show $x_*\in D(\Theta)$, we use $\xi_{n_k}\in \p \Theta(x_{n_k})$ to obtain
\begin{equation}\label{5.5.3.1}
\Theta(x_{n_k})\le \Theta(\hat{x}) +\l \xi_{n_k}, x_{n_k}-\hat{x}\r.
\end{equation}
In view of (\ref{12.15.11}) we have
\begin{equation*}
\Theta(x_{n_k})\le \Theta(\hat{x}) +\l \xi_{n_0}, x_{n_k}-\hat{x}\r + C_1 D_{\xi_{n_0}} \Theta(\hat{x}, x_{n_0}).
\end{equation*}
Since $x_{n_k}\rightarrow x_*$ as $k\rightarrow \infty$, by using the lower semi-continuity of $\Theta$ we obtain
\begin{equation}\label{5.5.2.1}
\fl \Theta(x_*)\le \liminf_{k\rightarrow\infty} \Theta(x_{n_k})\le
\Theta(\hat{x}) +\l \xi_{n_0}, x_*-\hat{x}\r + C_1 D_{\xi_{n_0}} \Theta(\hat{x}, x_{n_0})<\infty.
\end{equation}
This implies that $x_*\in D(\Theta)$.

In order to derive the convergence in Bregman distance, we first use (\ref{12.15.11}) to derive for $l<k$ that
\begin{eqnarray*}
\fl |\l\xi_{n_k}, x_{n_k}-x_*\r|
\le C_1 \left(D_{\xi_{n_l}} \Theta(x_*, x_{n_l}) -D_{\xi_{n_k}} \Theta(x_*, x_{n_k})\right)
+|\l \xi_{n_l}, x_{n_k}-x_*\r|.
\end{eqnarray*}
By taking $k\rightarrow \infty$ and using $x_{n_k}\rightarrow x_*$ we can derive that
\begin{eqnarray*}
\limsup_{k\rightarrow \infty} |\l\xi_{n_k}, x_{n_k}-x_*\r|
\le C_1 \left(D_{\xi_{n_l}} \Theta(x_*, x_{n_l}) - \varepsilon_0\right),
\end{eqnarray*}
where $\varepsilon_0:=\lim_{n\rightarrow \infty} D_{\xi_n}\Theta(x_*, x_n)$ which exists by the monotonicity of
$\{D_{\xi_n}\Theta(x_*, x_n)\}$. Since the above inequality holds for all $l$, by taking $l\rightarrow \infty$
we obtain
\begin{eqnarray}\label{5.5.6.1}
\limsup_{k\rightarrow \infty} |\l\xi_{n_k}, x_{n_k}-x_*\r|
\le C_1 \left(\varepsilon_0 - \varepsilon_0\right)=0.
\end{eqnarray}
Using (\ref{5.5.3.1}) with $\hat{x}$ replaced by $x_*$ we thus obtain
$\limsup_{k\rightarrow \infty} \Theta(x_{n_k})\le \Theta(x_*)$.
Combining this with (\ref{5.5.2.1}) we therefore obtain
\begin{equation*}
\lim_{k\rightarrow \infty} \Theta(x_{n_k})=\Theta(x_*).
\end{equation*}
This together with (\ref{5.5.6.1}) then implies that
\begin{equation*}
\lim_{k\rightarrow \infty} D_{\xi_{n_k}}\Theta(x_*, x_{n_k})=0.
\end{equation*}
Since $\{D_{\xi_n}\Theta(x_*, x_n)\}$ is monotonically decreasing, we can conclude that
\begin{equation*}
\lim_{n\rightarrow \infty} D_{\xi_n}\Theta(x_*, x_n)=0.
\end{equation*}

Finally we show that $x_*=x^\dag$. We use (\ref{5.5.3.1}) with $\hat{x}$ replaced by $x^\dag$ to obtain
\begin{equation}\label{5.5.7}
D_{\xi_0}\Theta(x_{n_k}, x_0)\le D_{\xi_0}\Theta(x^\dag, x_0)+\l \xi_{n_k}-\xi_0, x_{n_k}-x^\dag\r.
\end{equation}
By using (\ref{12.15.11}), for any $\varepsilon>0$  we can find $k_0$ such that
\begin{equation*}
\left|\l \xi_{n_k}-\xi_{n_{k_0}}, x_{n_k}-x^\dag\r\right| <\frac{\varepsilon}{2}, \qquad k\ge k_0.
\end{equation*}
We next consider $\l \xi_{n_{k_0}}-\xi_0, x_{n_k}-x^\dag\r$. Since
$\xi_{n+1}-\xi_n \in \overline{{\mathcal R}(F_0'(x^\dag)^*)}+\cdots + \overline{{\mathcal R}(F_{N-1}'(x^\dag)^*)}$,
we can find $v_{n,i}\in \Y_i^*$ and $\beta_{n,i} \in \X^*$ such that
\begin{equation*}
\fl \xi_{n+1}-\xi_n=\sum_{i=0}^{N-1} \left(F_i'(x^\dag)^* v_{n,i} +\beta_{n,i}\right) \quad \mbox{and} \quad \|\beta_{n,i}\|
\le \frac{\varepsilon}{3N B_1 n_{k_0}}, \quad 0\le n<n_{k_0},
\end{equation*}
where $B_1>0$ is a constant such that $\|x_n-x^\dag\|\le B_1$ for all $n$. Consequently
\begin{eqnarray*}
\fl \left|\l \xi_{n_{k_0}}-\xi_0, x_{n_k}-x^\dag\r \right|
&=\left|\sum_{n=0}^{n_{k_0}-1} \l \xi_{n+1}-\xi_n, x_{n_k}-x^\dag\r\right|\\
& =\left|\sum_{n=0}^{n_{k_0}-1} \sum_{i=0}^{N-1} \left[\l v_{n,i}, F_i'(x^\dag) (x_{n_k}-x^\dag)\r
+\l \beta_{n,i}, x_{n_k}-x^\dag\r \right]\right|\\
&\le \sum_{n=0}^{n_{k_0}-1} \sum_{i=0}^{N-1} \left(\|v_{n,i}\| \|F_i'(x^\dag) (x_{n_k}-x^\dag)\|
+\|\beta_{n,i}\| \|x_{n_k}-x^\dag\|\right)\\
&\le (1+\eta) \sum_{n=0}^{n_{k_0}-1} \sum_{i=0}^{N-1} \|v_{n,i}\| \|F_i(x_{n_k})-y_i\| +\frac{\varepsilon}{3}.
\end{eqnarray*}
Since $\|F_i(x_n)-y_i\|\rightarrow 0$ as $n\rightarrow \infty$, we can find $k_1\ge k_0$ such that
\begin{equation*}
|\l \xi_{n_{k_0}}-\xi_0, x_{n_k}-x^\dag\r |<\frac{\varepsilon}{2}, \qquad \forall k\ge k_1.
\end{equation*}
Therefore $|\l \xi_{n_k}-\xi_0, x_{n_k}-x^\dag\r|<\varepsilon$ for all $k\ge k_1$. Since $\varepsilon>0$ is arbitrary,
we obtain $\lim_{k\rightarrow \infty} \l \xi_{n_k}-\xi_0, x_{n_k}-x^\dag\r=0$. By taking $k\rightarrow \infty$ in
(\ref{5.5.7}) and using $\Theta(x_{n_k})\rightarrow \Theta(x_*)$ we obtain
\begin{equation*}
D_{\xi_0}\Theta(x_*, x_0)\le D_{\xi_0}\Theta(x^\dag, x_0).
\end{equation*}
According to the definition of $x^\dag$ we must have $D_{\xi_0}\Theta(x_*, x_0)=D_{\xi_0}\Theta(x^\dag, x_0)$.
A direct application of Lemma \ref{lem0} gives $x_*=x^\dag$. \hfill $\Box$

\begin{remark}
{\rm In the proof of Proposition \ref{general}, we obtain that $\lim_{k\rightarrow \infty} \Theta(x_{n_k})=\Theta(x_*)$. It is not clear if there holds
\begin{equation}\label{5.8.1}
\lim_{n\rightarrow \infty} \Theta(x_n) =\Theta(x_*)
\end{equation}
for the whole sequence $\{x_n\}$. Although the result of Proposition \ref{general} implies $x_n\rightarrow x_*$, we can not use it
to derive (\ref{5.8.1}) directly since $\Theta$ is not necessarily continuous at $x_*$.
}
\end{remark}

\section{Numerical examples}

In this section we will present some numerical simulations on Algorithm \ref{alg}. A key ingredient in this algorithm is the
resolution of the minimization problem
\begin{equation}\label{eq4.16.1}
x=\arg \min_{z\in \X} \left\{\Theta(z) -\l \xi, z\r\right\}
\end{equation}
for any given $\xi\in \X^*$. For some choices of $\Theta$, this minimization problem can be easily solved numerically.
In particular, when $\X=L^2(\Omega)$ and
\begin{equation}\label{eq:L1}
\Theta(x):= \frac{1}{2\beta} \int_\Omega |x(\omega)|^2 d\omega +\int_\Omega |x(\omega| d\omega
\end{equation}
with $\beta>0$, the minimizer of (\ref{eq4.16.1}) can be given explicitly by the soft thresholding
\begin{equation*}
x(\omega)=\left\{\begin{array}{lll}
\beta(\xi(\omega)-1) & \mbox{ if } \xi(\omega) >1,\\
0 & \mbox{ if } |\xi(\omega)|\le 1,\\
\beta(\xi(\omega) +1) & \mbox{ if } \xi(\omega) <-1.
\end{array}\right.
\end{equation*}
For the total variation like functional
\begin{equation}\label{eq:TV}
\Theta(x):= \frac{1}{2\beta} \int_\Omega |x(\omega)|^2 d\omega +\int_\Omega |Dx|
\end{equation}
with $\beta>0$ in $\X:=L^2(\Omega)$, although there is no explicit formula for the minimizer of (\ref{eq4.16.1}),
there are many numerical solvers developed in the literature; in our numerical simulations,
we will use the monotone version of the fast iterative shrinkage/thresholding algorithm (MFISTA) introduced
in \cite{BT2009}.

\subsection{Photoacoustic/thermoacoustic tomography}

Let $\Omega\subset {\mathbb R}^2$ be a bounded convex domain in ${\mathbb R}^2$ with smooth
boundary $\p \Omega$. We consider the problem of recovering a function $f : {\mathbb R}^2 \to {\mathbb R}$
supported over $\Omega$ from its means
\begin{equation*}
({\mathcal M} f) (x, r) :=\frac{1}{2\pi} \int_{{\mathbb S}^1} f (x + r\sigma) d\sigma
\end{equation*}
over the unit circle ${\mathbb S}^1$ with centers $x\in \p \Omega$ and radii $r > 0$. The resolution of this
problem is a key ingredient of many modern imaging techniques, including the photoacoustic tomography (PAT)
and thermoacoustic tomography (TAT).

PAT/TAT is a hybrid imaging technique based on the photoacoustic effect which refers to the generation of
acoustic waves by the absorption of electromagnetic (EM) energy coming from visible light, radio waves, or microwaves.
This technique combines the ultrasonic resolution with the high EM contrast to overcome the low contrast of
pure ultrasound diagnostics and weak resolution of pure optical imaging methodologies in soft tissues
(see \cite{KK2008,XW2006} and references therein). It provides diagnostics by determining the EM energy absorption
function $f(x)$ which strongly depends on the local biological properties of the cells. Upon absorption of
a short EM pulse, it results in thermoelastic expansion and thus in acoustic waves $p(x, t)$ that can be measured
by ultrasound receivers placed on a surface $S$ around the object. These measured data $p(x, t)$, $x\in S$ are then
used to determine $f(x)$. By assuming that the sound velocity is equal to the unit constant, the acoustic
wave $p(x,t)$ satisfies the Cauchy problem for the wave equation
\begin{eqnarray*}
\left\{\begin{array}{lll}
\p_t^2 p  -\Delta p =0 \qquad \mbox{ in } \,\,\,\, {\mathbb R}^3\times \{t\ge 0\},\\
p(x, 0)= f(x), \,\, \quad \p_t p(x,0)=0, \quad x\in {\mathbb R}^3.
\end{array}\right.
\end{eqnarray*}

The 3D TAT/PAT imaging can be reduced to 2D problems by introducing integrating line detectors \cite{BBGHP2007}.
Assuming the parallel detectors are perpendicular to the $x_2$-axis. Let the direction of detectors be $x_3$-axis,
the detectors measure $\bar p(x', t)$, $x':=(x_1, x_2)\in C$, for some detection curve $C$ in the $x_1x_2$-plane,
where $\bar p (x',t)$ denotes the projection of $p(x,t)$ in the $x_3$-direction
defined by
\begin{equation*}
\bar p(x',t) := \int_{-\infty}^\infty p(x', x_3, t) d x_3
\end{equation*}
and satisfies the 2D wave equation
\begin{eqnarray}\label{2Dwave}
\left\{\begin{array}{lll}
\p_t^2 \bar p -\Delta \bar p =0 \qquad \mbox{ in } \, \, {\mathbb R}^2 \times \{t\ge 0\},\\
\bar p(x', 0)=\bar f(x'), \,\,\quad \p_t \bar p(x', 0) =0, \quad x'\in {\mathbb R}^2,
\end{array}\right.
\end{eqnarray}
where $\bar f(x')=\int_{-\infty}^\infty f(x', x_3) d x_3$. Once $\bar f(x')$ can be determined from $\bar p(x',t)$, $x'\in C$,
we may rotate the detectors around the $x_2$-axis to obtain $\int_L f(x) d \ell$ for all lines $L$ perpendicular to $x_2$ axis
 which enable us to use the inversion
of Radon transform to reconstruct $f(x)$.

In order to determine $\bar f(x')$ from $\bar p(x',t)$, $(x',t)\in C\times [0, \infty)$, we use the representation formula for the solution
of (\ref{2Dwave}) to obtain (\cite{Sogge})
\begin{equation*}
\bar p(x', t) = \p_t \int_0^t \frac{r ({\mathcal M} \bar f)(x',r)}{\sqrt{t^2-r^2}} d r.
\end{equation*}
By using the solution formula of Abel integral equation, it follows (\cite{BBGHP2007})
\begin{equation*}
({\mathcal M} \bar f)(x', r) = \frac{2}{\pi} \int_0^r \frac{\bar p(x',t)}{\sqrt{r^2-t^2}} dt, \quad x'\in C \mbox{ and } r\ge 0.
\end{equation*}
Therefore we need to determine $\bar f(x')$ from its means $({\mathcal M}\bar f)(x', r)$, $x'\in C$ and $r\ge 0$.

\begin{figure}[htp]
  \begin{center}
    \includegraphics[width = 1\textwidth, height= 3in]{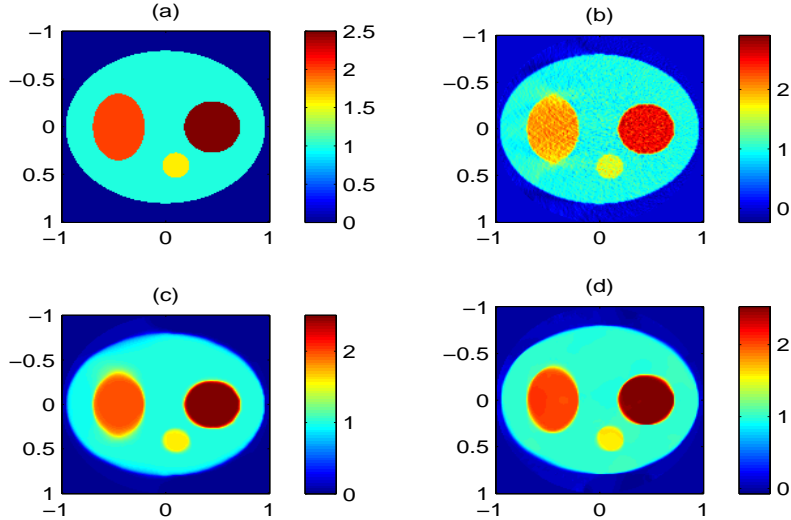}
  \end{center}
  \caption{Numerical results with $N=80$ measurements and $\d=0.01$: (a) exact solution;
  (b) $\Theta(f) =\|f\|_{L^2}^2$; (c) and (d) $\Theta(f) = \frac{1}{2\beta} \|f\|_{L^2}^2 +\int_\Omega |Df|$
  with $\beta=1$ and $\beta=10$ respectively}\label{fig1}
\end{figure}

In numerical simulations, we reconstruct a function $f$ supported on the disk $B_R$ of radius $R=0.96$ centered at the origin
from its means $({\mathcal M} f)({\bf x}_j, r)$, $r\ge 0$, measured at $N$ points ${\bf x}_j =R(\sin(j\pi/N), \cos(j\pi/N))$, $j=0, \cdots, N-1$
uniformly distributed on the semicircle $S^{+}:= \{{\bf x}\in\partial B_R: x_1\geq0\}$. This is equivalent to solving the system
\begin{equation*}
 {\mathcal M}_i f = g_i, \quad i=0, \cdots, N-1,
\end{equation*}
where
\begin{equation*}
({\mathcal M}_i f)(r):= ({\mathcal M} f)({\bf x}_i, r)=\frac{1}{2\pi} \int_{{\mathbb S}^1} f({\bf x}_i+ r \sigma) d \sigma.
\end{equation*}
It is easy to check \cite{HKLS2007} that the operators ${\mathcal M}_i$ can be continuously extended to
\begin{equation*}
{\mathcal M}_i: L^2(B_R)\to L^2([0,2R],rdr)
\end{equation*}
with $\|{\mathcal M}_i\|\leq 2\sqrt{\pi}$ and the adjoint ${\mathcal M}_i^*: L^2([0,2R],rdr)\to L^2(B_R)$ is given by
 \begin{equation*}
({\mathcal  M}_i^*g)({\bf x})= 2 g(|{\bf x}_i-{\bf x}|).
 \end{equation*}

In Figure \ref{fig1} we report the numerical results with $N = 80$ measurements. In order to approximate functions, we take $160\times 160$ grid
points uniformly distributed on the square $[-1,1]\times [-1,1]$. The exact piecewise constant phantom $f^\dag$ is shown in (a).
For the simulations, we add 2\% uniformly distributed noise to ${\mathcal M}_i f^\dag$ to produce the noisy data which are then used
to reconstruct $f^\dag$. Figure \ref{fig1} (b)--(d) give the reconstruction result by Algorithm \ref{alg} with initial guess $f_0=\xi_0=0$ and
different choices of the uniformly convex functionals $\Theta$, we take $r=2$, $\tau=1.2$ and $\mu_0=(1-1/\tau)/(\beta\sqrt{\pi})$ when $\Theta$
given by (\ref{eq:L1}) and (\ref{eq:TV}) are used. Figure \ref{fig1} (b) presents the reconstruction
result by Algorithm \ref{alg} with $\Theta(f) =\int_{B_R} |f({\bf x})|^2 d {\bf x}$. Figure \ref{fig1} (c) and (d) report the reconstruction results
by Algorithm \ref{alg} with $\Theta(f)$ defined by (\ref{eq:TV}) with $\beta =1$ and $\beta =10$ respectively.
The results in (c) and (d) significantly improve the one in (b) by efficiently removing the artefacts due to the notorious
oscillatory effect and indicate that the results are robust with respect to $\beta$.

\begin{figure}[htp]
  \begin{center}
    \includegraphics[width = 1\textwidth, height= 1.3in]{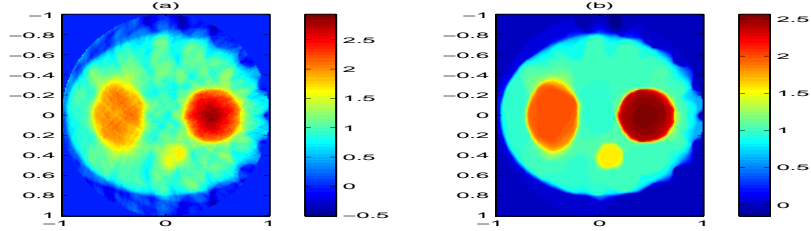}
  \end{center}
  \caption{Numerical results with $N=10$ measurements and $\d=0.01$: (a) $\Theta(f) =\|f\|_{L^2}^2$;
  (b) $\Theta(f) = \frac{1}{2} \|f\|_{L^2}^2 +\int_\Omega |Df|$}\label{fig2}
\end{figure}

In Figure \ref{fig2} we report the computation results by reducing the number of measurements to $N=10$. The reconstruction ability of the method using
$\Theta(f)=\|f\|_{L^2}^2$ becomes worse, the method with $\Theta(f)=\frac{1}{2} \|f\|_{L^2}^2 +\int_\Omega |Df|$,
however, still has good reconstruction. This reflects the philosophy in compressed sensing: it is possible to reconstruct an image
from very few number of measurements by using $L^1$ penalty term if it is sparse under suitable transformation.

\subsection{Parameter identification}

We next consider the identification of the parameter $c$ in the boundary value problem
\begin{eqnarray}\label{PDE}
\left\{
  \begin{array}{ll}
    -\triangle u + cu = f \qquad  \mbox{in } \Omega, \\
    u = g  \qquad \mbox{on } \partial\Omega,
  \end{array}
\right.
\end{eqnarray}
from an $L^2(\Omega)$-measurement of the state $u$, where $\Omega\subset\mathbb{R}^d$ with $d\leq 3$ is
a bounded domain with Lipschitz boundary, $f\in L^2(\Omega)$ and $g\in H^{3/2}(\Omega)$. We assume that
the exact solution $c^{\dag}$ is in
$L^2(\Omega)$. This problem reduces to solving $F(c) =u$, i.e. (\ref{sys}) with $N=1$, if we define the
nonlinear operator $F: L^2(\Omega)\rightarrow  L^2(\Omega)$ by
\begin{equation}\label{nonlinear}
F(c): = u(c)\,,
\end{equation}
where $u(c)\in H^2(\Omega)\subset  L^2(\Omega)$ is the unique solution of (\ref{PDE}). This operator $F$
is well defined on
\begin{equation*}
D(F): = \left\{ c\in  L^2(\Omega) : \|c-\hat{c}\|_{ L^2(\Omega)}\leq \gamma_0 \mbox{ for some }
\hat{c}\geq0,\ \textrm{a.e.}\right\}
\end{equation*}
for some positive constant $\gamma_0>0$. It is known that $F$ is Fr\'{e}chet differentiable; the Fr\'{e}chet derivative
of $F$ and its adjoint are given by
\begin{equation}\label{Frechet}
F'(c)h  = -A(c)^{-1}(hF(c)) \quad \mbox{and} \quad F'(c)^* w  = -u(c) A(c)^{-1}w
\end{equation}
for $h, w\in L^2(\Omega)$, where $A(c): H^2\cap H_0^1\rightarrow  L^2$ is defined by $A(c) u = -\triangle u +cu$
which is an isomorphism uniformly in ball $B_{\rho}(c^{\dag}) \cap D(F)$ for small $\rho>0$.
Moreover, Assumption \ref{A2} holds for small $\rho>0$ (see \cite{HNS95}).

\begin{figure}[htp]
  \begin{center}
    \includegraphics[width = 1\textwidth, height= 3.5in]{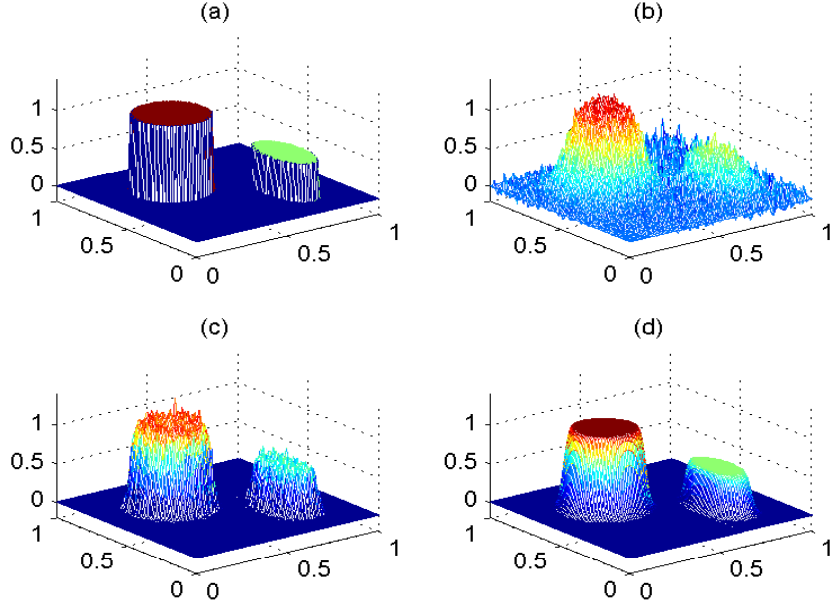}
  \end{center}
  \caption{ numerical results with $\d=0.5\times 10^{-4}$ and $\tau=1.1$: (a) exact solution;
  (b) $\Theta(f) =\|f\|_{L^2}^2$; (c) $\Theta(c) = \frac{1}{2}\|c\|_{L^2}^2 +\|c\|_{L^1}$;
  (d) $\Theta(f) = \frac{1}{2} \|f\|_{L^2}^2 +\int_\Omega |Df|$}\label{fig3}
\end{figure}

In our numerical simulation, we consider the two dimensional problem with $\Omega = [0,1]\times[0,1]$ and
\begin{eqnarray*}
c^{\dag}(x,y) = \left\{
  \begin{array}{ll}
    1, \qquad & \hbox{if } (x-0.65)^2 + (y-0.36)^2 \le 0.18^2, \\
    0.5, & \hbox{if } (x-0.35)^2 + 4 (y-0.75)^2 \le 0.2^2, \\
    0, & \hbox{elsewhere}.
  \end{array}
\right.
\end{eqnarray*}
We assume $u(c^{\dag}) = x+y$ and add random noise to produce $u^{\delta}$ satisfying
$\|u^{\delta} - u(c^{\dag})\|_{L^2(\Omega)} = \delta$ with $\d=0.5\times 10^{-4}$. In order to reconstruct $c^{\dag}$,
we use Algorithm \ref{alg} with $r=2$, $N=1$ and $\tau=1.1$ and take $c_0=\xi_0=0$ as initial guess; we take the step length
$\mu_{n,i}^\d$ according to Remark \ref{rem:6.21.1} with
$$
\tilde{\mu}_{n,i}^\d = \min\left\{\frac{\mu_0 \|F_i(c_{n,i}^\d)-y_i^\d\|_{L^2}^2}{\|F_i'(c_{n,i}^\d)^*(F_i(c_{n,i}^\d)-y_i^\d)\|_{L^2}^2}, \mu_1\right\}
$$
with $\mu_0=(1-1/\tau)/\beta$ and $\mu_1=4000$ when $\Theta$ given by (\ref{eq:L1}) and (\ref{eq:TV}) are used.
The partial differential equations involved are solved approximately by a finite difference method by dividing $\Omega$
into $100\times 100$ small squares of equal size.

We report the numerical results in Figure \ref{fig3}. In (a) we plot the exact solution $c^{\dag}(x,y)$. In (b) we plot the result
of Algorithm \ref{alg} with $\Theta(c) = \|c\|^2_{L^2}$; although the reconstruction tells something on the sought solution,
it does not tell more information such as sparsity, discontinuities and constancy since the result is too oscillatory. In (c) we
report the result of Algorithm \ref{alg} with $\Theta$ given by (\ref{eq:L1}) with $\beta=1$. It is clear that the sparsity
of the sought solution is significantly reconstructed. The reconstruction result, however, is still oscillatory
on the nonzero parts which is typical for this choice of $\Theta$. In (d) we report the result of Algorithm \ref{alg}
with $\Theta(c)$ given by (\ref{eq:TV}) with $\beta =1$. The reconstruction is rather satisfactory and
the notorious oscillatory effect is efficiently removed.

\subsection{Schlieren imaging}

Consider the problem of reconstructing a function $f$ supported on a bounded domain $D\subset {\mathbb R}^2$
from
\begin{equation*}
{\mathcal I} f(s,\sigma) =\left(\int_{\mathbb R} f(s \sigma +r \sigma^\perp) d r \right)^2, \qquad (s, \sigma) \in {\mathbb R}\times {\mathbb S}^1.
\end{equation*}
This problem arises from determining the 3D pressure fields on cross-sections of a water tank generated by an ultrasound transducer from Schlieren data.
The data are collected with a Schlieren optical system based on Raman scattering. The Schlieren optical system outputs the intensity of light through the tank
which is proportional to the square of the line integral of the pressure along the light path \cite{HZ91}.

In our numerical simulations we reconstruct a function $f$ support on
$
D=[-1,1]\times [-1,1]
$
from $N=100$ recording angles $\sigma_i \in {\mathbb S}^1$, $i = 0, \cdots, N-1$ uniformly distributed on the semicircle.
It reduces to solving the system (\ref{sys}) by introducing
\begin{equation*}
F_i(f)(s): = {\mathcal I} f(s,\sigma_i):= \left({\mathcal R}_i f(s)\right)^2,
\end{equation*}
where ${\mathcal R}_i f(s) =\int_{\mathbb R} f(s \sigma_i +r \sigma_i^\perp) d r$ is the Radon transform.
It is easy to show (\cite{HKLS2007}) that each $F_i: H_0^1(D) \to L^2(I)$ with $I=[-\sqrt{2},\sqrt{2}]$ is Fr\'{e}chet differentiable with
\begin{equation*}
F_i'(f) h ={\mathcal R}_if\cdot {\mathcal R}_i h, \quad \forall h\in H_0^1(D)
\end{equation*}
and the adjoint $F_i'(f)^*: L^2(I) \to H_0^1(D)$ is given by
\begin{equation*}
F_i'(f)^* g = (I-\Delta)^{-1} (2 {\mathcal R}_i^\# (g {\mathcal R}_i f)), \quad \forall g\in L^2(I),
\end{equation*}
where ${\mathcal R}_i^\#: L^2(I)\to L^2(D)$ is the adjoint of ${\mathcal R}_i$ given by  $({\mathcal R}_i^\# g)(x)= g(\l x, \sigma_i\r)$.

\begin{figure}[htp]
  \begin{center}
    \includegraphics[width = 1\textwidth, height= 3in]{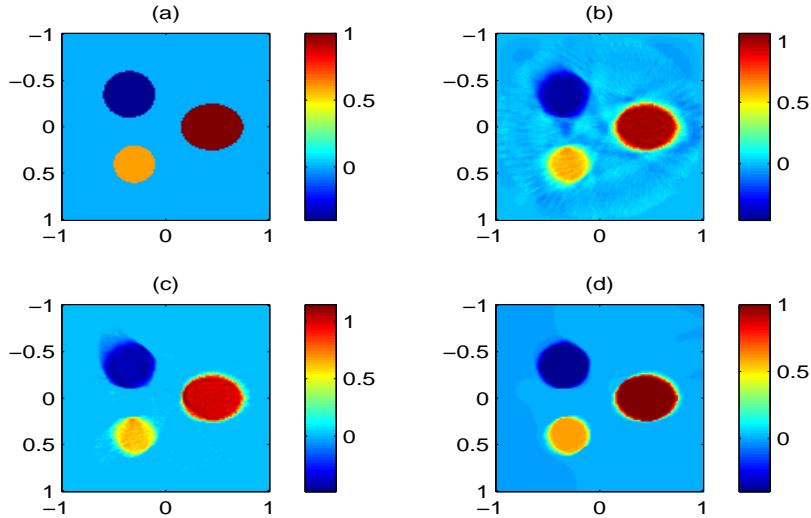}
  \end{center}
  \caption{ Numerical results with $N=100$ measurements, $\d=0.002$ and $\tau=1.5$: (a) exact solution;
    (b) $\Theta(f)=\|f\|_{L^2}^2$; (c) $\Theta(f) = \frac{1}{2}\|f\|_{L^2}^2 +\|f\|_{L^1}$;
    (d) $\Theta(f)= \frac{1}{2} \|f\|_{L^2}^2 +\int_\Omega |Df|$}\label{fig4}
\end{figure}

In Figure \ref{fig4} we report the numerical results with $N = 100$ measurements with the exact solution $f^\dag$ shown in (a).
For the simulations, we use noisy data with noise level $\d=0.002$ to reconstruct $f^\dag$; in order to approximate functions, we take
$120\times 120$ grid points uniformly distributed over $[-1,1]\times[-1,1]$. Figure \ref{fig4} (b)--(d) give the
reconstruction results by Algorithm \ref{alg} implemented with $r=2$ and $\tau =1.5$, initial guess $\xi_0=0.01$ and different choices of $\Theta$;
the step length $\mu_{n,i}^\d$ is chosen according to Remark \ref{rem:6.21.1} with $\mu_0=(1-1/\tau)/\beta$ and $\mu_1=1000$. Figure \ref{fig4} (b)
clearly contains many artefacts. Figure \ref{fig4} (c) removes almost all artefacts since it uses the $L^1$ like penalty functional,
the reconstruction on nonzero profile, however, turns out to be unsatisfactory. Because of the use of total variation like penalty functional,
Figure \ref{fig4} (d) efficiently removes the artefacts and reconstructs the profile very well.

%\begin{figure}[ht]
%\centering
%  \includegraphics[width=5in, height=2.5in]{L1_1d.eps}
%  \caption{Reconstruction results for Example \ref{ex1} with $n_\d$ being the number
%  of iterations. (a) $n_\d=16$; (b) $n_\d=16$; (c) $n_\d=24$.}\label{fig1}
%\end{figure}

\section*{Acknowledgement}

The work of Q Jin is partially supported by the DECRA grant DE120101707 of Australian Research Council.

%%%%%%%%%%%%%%%%%%%%%%%%%%%%%%%%%%%%%%%%%%%%%%%%%%%%%%%%%%%%%%%%%%%%%%%%%%%
%

\end{document}